\documentclass[10pt,oneside,a4paper]{amsart}
\usepackage{amsmath,amsfonts,amsthm, amssymb}
\usepackage{array}
\usepackage[all,textures]{xy}

\theoremstyle{plain}
\newtheorem{theorem}{Theorem}[section]
\newtheorem{proposition}[theorem]{Proposition}
\newtheorem{lemma}[theorem]{Lemma}
\newtheorem{corollary}[theorem]{Corollary}

\theoremstyle{definition}
\newtheorem{definition}[theorem]{Definition}

\theoremstyle{remark}

\newtheorem{remark}[theorem]{Remark}

\newcommand{\ovl}{\overline}
\newcommand{\Kern}{\mathrm{Ker}}

\newcommand{\D}{\mathrm{D}}

\renewcommand{\lim}{\mathrm{lim}}

\newcommand{\Tor}{\mathrm{Tor}}
\newcommand{\Ext}{\mathrm{Ext}}
\newcommand{\Hom}{\mathrm{Hom}}

\newcommand{\Sk}{\mathrm{Sk}}

\newcommand{\Z}{\mathbb{Z}}

\newcommand{\N}{\mathbb{N}}

\newcommand{\CCC}{\mathfrak{c}}

\newcommand{\Set}{\ensuremath{\mathsf{Set}} }

\newcommand{\Ab}{\ensuremath{\mathsf{Ab}} }
\newcommand{\Mod}{\ensuremath{\mathsf{Mod}} }

\newcommand{\mmod}{\ensuremath{\mathsf{mod}} }

\newcommand{\Gd}{\ensuremath{\mathsf{Gd}} }

\newcommand{\Ind}{\ensuremath{\mathsf{Ind}}}
\newcommand{\Fun}{\ensuremath{\mathsf{Fun}}}
\newcommand{\Add}{\ensuremath{\mathsf{Add}}}
\newcommand{\Inj}{\ensuremath{\mathsf{Inj}}}
\newcommand{\Cof}{\ensuremath{\mathsf{Cof}}}

\newcommand{\lra}{\longrightarrow}

\newcommand{\ccc}{\ensuremath{\mathcal{C}}}
\newcommand{\ddd}{\ensuremath{\mathcal{D}}}

\newcommand{\LLL}{\ensuremath{\mathcal{L}}}

\CompileMatrices
\SelectTips{cm}{11}

\title{Obstruction theory for objects in abelian and derived
categories}
\author{Wendy T. Lowen$^{\ast}$}
\address{Departement DWIS\\ Vrije Universiteit Brussel\\ Pleinlaan
2\\1050 Brussel\\ Belgium}
\email[Wendy T. Lowen]{wlowen@vub.ac.be}
\thanks{$^{\ast}$aspirant at the FWO}

\begin{document}
\begin{abstract}
  In this paper we develop the obstruction theory for lifting
  complexes, up to quasi-isomorphism, to derived categories of flat
  nilpotent deformations of abelian categories.  As a particular case
  we also
  obtain the corresponding obstruction theory for lifting of objects
  in terms of Yoneda $\Ext$-groups.
  In appendix we prove the existence of miniversal derived deformations of
  complexes.
\end{abstract}
\maketitle
\section{Introduction}
Complete families of non-commutative deformations of projective
planes, quadrics and more generally Hirzebruch surfaces where
constructed in \cite{Bondal,VdB26,VdB28} using adhoc deformation theoretic
arguments. In order to provide a firmer foundation for these
constructions we developed in \cite{defab,hoch} a deformation theory for
abelian categories which generalizes the deformation theory of (module
categories over) algebras.

The arguments in \cite{Bondal,VdB26,VdB28} are based on the intuition that
exceptional objects \cite{Bondal2}  should lift to
any deformation.  In the current paper we will justify this assumption
by developing an obstruction theory for the lifting of objects (and
complexes) to deformations of an abelian category.

Let us first summarize the deformation theory of abelian categories.
Assume that $R \lra R_0$ is a surjective ringmap with nilpotent kernel
between coherent, commutative rings\footnote{In applications $R$ and
  $R_0$ will probably be artinian local rings but the added generality
  we allow incurs very little cost}. A deformation of an $R_0$-linear
abelian category $\ccc_0$ along $R \lra R_0$ is an $R$-linear functor
$\ccc_0 \lra \ccc$ inducing an equivalence $\ccc_0 \cong \ccc_{{R_0}}$
where $\ccc_{{R_0}} \subset \ccc$ is the full subcategory of
${{R_0}}$-objects, i.e.\ objects with an ${{R_0}}$-structure
compatible with the $R$-structure \cite[Def.5.2, \S4]{defab}. In
general such deformations can be very wild but we show in loc.\ cit.
that by restricting to (appropriately defined) flat 
deformations the theory becomes controllable.

The definition of flatness for an abelian $R$-linear category is
somewhat involved \cite[Def.3.2]{defab} but for a category with enough
injectives it amounts to requiring that injectives are $R$-coflat, i.e.\
$R$-flat in the dual category \cite{AZ2}.   An $R_0$-algebra $A_0$ is
flat if and only if $\Mod(A_0)$ is flat, and flat $R$-deformations of
$\Mod(A_0)$ correspond precisely to flat $R$-deformations of $A_0$ \cite{defab}.

\def\RHom{\operatorname{RHom}}
\def\Lotimes{\operatorname{\overset{L}{\otimes}}}

In this paper we will study the problem of lifting objects along the
functor $\Hom_R({{R_0}},-):\ccc\lra\ccc_0$ for a deformation $\ccc_0 \lra \ccc$ and
similarly the problem of lifting objects in the correspoding derived
categories along the functor $\RHom_R({{R_0}},-)$. By dualizing one
obtains lifting properties for the (perhaps more familiar) functors
${{R_0}}\otimes_R-$ and ${{R_0}}\Lotimes_R-$. We leave the explicit
formulations of these dual versions to the reader.  There is a
parallel obstruction theory for lifting maps which is contained in the
body of the paper, but which for brevity we will not formulate in this
introduction.

Consider
surjective ringmaps between coherent, commutative rings $$\bar{R} \lra
R \lra R_0$$
with $\Kern(\bar{R} \lra {{R_0}}) = I$, $\Kern(\bar{R} \lra
R) = J$ and $IJ = 0$.  Consider flat abelian deformations $\bar{\ccc}
\longleftarrow \ccc \longleftarrow \ccc_0$ along these ring maps along
with their adjoints
$$\Hom_{\bar{R}}(R,-): \bar{\ccc} \lra \ccc
\,\,\,\,\text{and}\,\,\,\, \Hom_R({{R_0}},-): \ccc \lra \ccc_0
$$
For a functor $F$ and an object $C$ in the codomain of $F$,
$L_F(C)$ denotes the natural groupoid of lifts of $C$ along $F$
(Definition \ref{deflift1}).

We prove
the following obstruction theory for lifting coflat objects along the
restricted functor
$$\Hom_{\bar{R}}(R,-): \Cof(\bar{\ccc}) \lra
\Cof(\ccc)$$
where $\Cof(-)$ denotes the full subcategory of coflat objects. 
{\def\thetheorem{A}
\begin{theorem}\label{Klift5intro}
Consider a lift
$C$ of $C_0$ along $\Hom_{{R}}({{R_0}},-)$.
\begin{enumerate}
\item There is an obstruction $$o(C) \in 
\Ext^2_{\ccc_0}(R\Hom_{{R_0}}(J,C_0),C_0)$$
with $$o(C) = 0 \iff L_{\Hom_{\bar{R}}(R,-)}(C) \neq \varnothing.$$
\item If $o(C) = 0$, then $\Sk(L_{\Hom_{\bar{R}}(R,-)}(C))$ is affine
over $$\Ext^1_{\ccc_0}(R\Hom_{{R_0}}(J,C_0),C_0).$$
\end{enumerate}
\end{theorem}
}
The previous result generalizes the classical obstuction theory for lifting
along
$\Hom_{\bar{R}}(R,-):\Mod(\bar{A}) \lra \Mod(R \otimes_{\bar{R}} \bar{A})$ for
an
$\bar{R}$-algebra $\bar{A}$.
\cite{Laudal}. Note that as expected when ${{R_0}}$ is a field, we obtain
obstructions purely in terms of the Yoneda $\Ext$-groups
$\Ext^i_{\ccc_0}(C_0,C_0)$.

Theorem $\ref{Klift5intro}$ is closely related to
our main Theorem
\ref{Klift4intro} below (which is contained in Theorem \ref{corcor} in the
body of the paper).  Theorem \ref{Klift4intro} gives the
obstruction theory for lifting along the restricted derived functor
$$R\Hom_{\bar{R}}(R,-): D^b_{{\mathrm {fcd}}}(\bar{\ccc}) \lra
D^b_{{\mathrm {fcd}}}(\ccc).$$
Here ``$\mathrm {fcd}$'' means that we restrict to objects of finite 
coflat dimension
(Definition \ref{defcd}). The dual of this condition is finite $\Tor$-dimension ($\mathrm {ftd}$), as considered for example in
\cite{De}.
{\def\thetheorem{B}
\begin{theorem}\label{Klift4intro}
Consider a lift $C^{\cdot}$ of $C^{\cdot}_0$ along $R\Hom_{{R}}({{R_0}},-)$.
\begin{enumerate}
\item There is an obstruction $$o(C^{\cdot}) \in
\Ext^2_{\ccc_0}(R\Hom_{{R_0}}(J,C^{\cdot}_0),C^{\cdot}_0)$$ with 
$$o(C^{\cdot}) = 0
\iff L_{R\Hom_{\bar{R}}(R,-)}(C^{\cdot}) \neq \varnothing.$$
\item If $o(C^{\cdot}) = 0$, then
$\Sk(L_{R\Hom_{\bar{R}}(R,-)}(C^{\cdot}))$ is affine over
$$\Ext^1_{\ccc_0}(R\Hom_{{R_0}}(J,C^{\cdot}_0),C^{\cdot}_0).$$
\end{enumerate}
\end{theorem}
}
Our approach for proving Theorem \ref{Klift4intro} is to look at the functor
\begin{equation}\label{label}
K(\Hom_{\bar{R}}(R,-)): K(\Inj(\bar{\ccc})) \lra K(\Inj(\ccc)
\end{equation}
between homotopy categories for a deformation with enough injectives,
which leads to the problem of
naively deforming differentials and cochain maps to fixed graded lifts
of complexes. A
detailed obstruction theory for this problem (Theorem
\ref{groot2}) is worked out in section \S\ref{parliftdifcoc} for a full
additive functor $F: \bar{\CCC} \lra \CCC$ with $(\Kern(F))^2 = 0$
(\S\ref{parobs} (\ref{kernil})) between additive categories.
For such a functor, we prove a ``crude lifting lemma'' inspired by the Crude
Perturbation Lemma in \cite{markl}, which implies that every lift of a complex
along
\begin{equation}\label{kc}
K(F): K(\bar{\CCC}) \lra K(\CCC)
\end{equation}
is homotopy equivalent to a lift of its differential to a fixed graded
lift.  If $F$ is essentially surjective, this leads to the obstruction
theory for such along $K(F)$ (Theorem \ref{Klift}).  Our main example
of a functor $F$ with the indicated properties is a linear deformation
along $\bar{R} \lra R$, i.e a $\bar{R}$-linear functor $\bar{\CCC}
\lra \CCC$ inducing $R \otimes_{\bar{R}} \bar{\CCC} \cong \CCC$. Here
$R \otimes_{\bar{R}} \bar{\CCC}$ is obtained from $\bar{\CCC}$ by
tensoring the hom-sets with $R$.
%A linear deformation $\bar{\CCC} \lra \CCC$ between small categories induces
%an abelian deformation $\Mod(\CCC) \lra \Mod(\bar{\CCC})$, linking the
%theories of linear and abelian deformations (\cite[Prop 5.4]{defab}).

Consider
flat linear deformations
$$\xymatrix{{\bar{\CCC}} \ar[r]_F & {\CCC} \ar[r]_{(-)_0} & {\CCC_0}}$$ along
$\bar{R} \lra R \lra {{R_0}}$ (here flat means that the hom-sets are 
flat modules).
In Theorem \ref{Klift2}, we show that the lifts of a complex $C^{\cdot} \in
K(\CCC)$ along $K(F)$ are governed by the complex
\begin{equation}\label{xxx}
J \otimes_{{R_0}} \Hom_{\CCC_0}(C_0^{\cdot},C_0^{\cdot}).
\end{equation}
Since for a flat abelian deformation $\bar{\ccc} \longleftarrow \ccc$ with
enough injectives,
$\Hom_{\bar{R}}(R,-): \Inj(\bar{\ccc}) \lra \Inj(\ccc)$ defines a flat linear
deformation (Proposition \ref{nilpinjlift}), the complex (\ref{xxx}) for
(\ref{label}) translates into
$$R\Hom_{\ccc_0}(R\Hom_{{R_0}}(J,C_0),C_0),$$
which is the complex behind Theorem \ref{Klift4intro}.

For completeness, we prove the existence of miniversal homotopy and derived
deformations of complexes (when we consider trivial linear or abelian
deformations of categories) in Appendix, using Schlessingers conditions
\cite{schlessinger}.

\medskip

To the best of the author's knowledge Theorems A and B (and their
generalizations to maps stated below) have not been formulated in
the current generality before. However some particular cases are certainly
known. For Theorem A we have already mentioned module categories
\cite{Laudal}. The case of coherent sheaves over algebraic varieties
is also standard (see for example \cite{Vistoli}). First order
deformations of an object (for a trivial deformation of an abelian
category) were classified in \cite{AZ2}.  Theorem B was proved
by Inaba for the derived category of coherent sheaves over a
projective variety \cite{Inaba}. Related results for the derived
category of a profinite group are stated in \cite{defder}.

\medskip

The author wishes to thank Michel Van den Bergh for suggesting the use
of injective resolutions and for several interesting discussions.

\section{Notations and preliminaries on cochain complexes}\label{not}

Let $\ccc$ be a linear category, i.e. a category enriched over some module
category.  We have the
\emph{graded category}
$G(\ccc) =
\Fun(\Z,\ccc)$ whose objects are denoted by ${C}^{\cdot}$ and are called
\emph{graded objects}. For $C^{\cdot} \in G(\ccc)$, $C^{\cdot}[n]$ denotes the
shifted graded object with $C[n]^i = C^{i+n}$. A morphism $C^{\cdot} \lra
D^{\cdot}[n]$ is called a \emph{graded map of degree $n$ from $C^{\cdot}$
to $D^{\cdot}$}. The composition of a graded map of degree $n$ with a
graded map of degee $m$ is a graded map of degree $n+m$. A graded
map of degree
$1$ from
$C^{\cdot}$ to $C^{\cdot}$ will be called a \emph{pre-differential on
$C^{\cdot}$}. A \emph{differential on $C^{\cdot}$} is a pre-differential $d$
with $d^2 = 0$.
For graded objects $C^{\cdot}, D^{\cdot} \in G(\ccc)$, we define a graded
abelian group $\Hom^{\cdot}(C^{\cdot}, D^{\cdot})$ by
$\Hom^n(C^{\cdot}, D^{\cdot}) = G(\ccc)(C^{\cdot}, D^{\cdot}[n])$.
For pre-differentials $d_C$ on $C^{\cdot}$ and $d_D$ on $D^{\cdot}$, we define
the pre-differential $\delta = \delta_{d_C,d_D}$ on  $\Hom^{\cdot}(C^{\cdot},
D^{\cdot})$
by $$\delta^n(f) = d_Df - (-1)^nfd_C.$$
If $d_C = d_D$,
$\delta$ turns $\Hom^{\cdot}(C^{\cdot}, C^{\cdot})$ into a cDG-algebra
\cite{Schwarz}.
If $d_C$ and $d_D$ are differentials, then so is $\delta$. In this case, if
$d_C = d_D$, $\Hom^{\cdot}(C^{\cdot}, C^{\cdot})$ becomes a DG-algebra.
A \emph{pre-complex} $(C^{\cdot},d)$ is a graded object $C^{\cdot}$ endowed
with a pre-differential $d$, if $d$ is a differential then $(C^{\cdot},d)$ is
called a \emph{(cochain) complex}. A \emph{cochain map of degree
$n$} between pre-complexes
$(C^{\cdot},d_C) \lra (D^{\cdot},d_D)$ is a graded map $f \in
\Hom^n(C^{\cdot}, D^{\cdot})$ with
$\delta(f) = 0$.
For graded maps $f, g \in \Hom^n(C^{\cdot},
D^{\cdot})$, a \emph{homotopy} $H: f \lra g$ is a graded map $H \in
\Hom^{n-1}(C^{\cdot}, D^{\cdot})$ with $\delta(H) = g - f$.
Pre-complexes, cochain maps and homotopies constitute a bicategory $P(\ccc)$
in which the complexes form a full bisubcategory $C(\ccc)$.
The \emph{homotopy category} $K(\ccc)$ is obtained from $C(\ccc)$ by
considering cochain maps up to homotopy. Restricting to bounded below
complexes yields the category $K^+(\ccc)$. If $\ccc$ is an \emph{abelian}
category, there is a functor
$C(\ccc)
\lra G(\ccc): C^{\cdot}
\longmapsto H^{\cdot}C$ mapping a cochain complex to its graded homology
object. Cochain maps which are mapped onto isomorphisms by this functor are
called
\emph{quasi-isomorphisms}. The derived category $D(\ccc)$ is obtained from
$C(\ccc)$ by formally inverting all quasi-isomorphisms. Restricting to bounded
below or bounded complexes yields the derived categories
$D^+(\ccc)$ and $D^b(\ccc)$ respectively.

\section{Lifting differentials and cochain maps}\label{parliftdifcoc}

In \S \ref{parobs}, we develop the obstruction theory for naively lifting
differentials and cochain maps along a
suitable additive functor, relative to fixed graded lifts
(Corollary \ref{obscoc} and
Theorem \ref{groot2}). \S \ref{parcrude} contains some comparison results
for the obstructions defined in \S \ref{parobs}, which enable us to prove
a ``crude homological lifting lemma'' (Corollary \ref{crude}) which refers to
the Crude Perturbation Lemma in \cite{markl}. Since the Crude Perturbation
Lemma does not immediately apply, we give a proof of Corollary \ref{crude} in
this paper. However, we believe a generalization of \cite{markl} to
perturbations of ``complexes-modulo-a-subcategory''would also capture Corollary
\ref{crude}. We start with introducing some terminology.

\subsection{Some lift groupoids}\label{parliftobmap}\label{pardefcoc}
In this section we define the various lift groupoids we will use throughout
this paper.
Let $F: \ovl{\ccc} \lra \ccc$ be an arbitrary functor.
\begin{definition}\label{deflift1}
\begin{enumerate}
\item For an object $C \in
\ccc$, a \emph{lift of $C$ along $F$} is an object $\ovl{C} \in 
\ovl{\ccc}$ together
with an isomorphism
$c: C
\cong F(\ovl{C})$. A lift $(\ovl{C},c)$ of $C$ will often be denoted simply by
$\ovl{C}$ or $c$. If $F:
\ovl{\ccc}
\lra
\ccc$ is right adjoint to a functor
$G:
\ccc
\lra
\ovl{\ccc}$, a lift of $C$ along $F$ can be represented by a map $G(C) \lra
\ovl{C}$.
\item For a map $f: C
\lra C'$ in $\ccc$ and lifts $c:C \cong F(\ovl{C})$ and $c':C' \cong
F(\ovl{C'})$ along $F$ of
$C$ and
$C'$ respectively, a \emph{lift of $f$ (along $F$) relative to $c,c'$} is a
map $\ovl{f}:
\ovl{C}
\lra
\ovl{C'}$ with
$F(\ovl{f})c = c'f$. The set of all lifts of $f$ along $F$ relative to
$c,c'$ will be denoted by $$L_F(f|c,c').$$
\item For $C \in \ccc$, we consider the following groupoid $$L_F(C):$$
\begin{enumerate}
\item[0.] Objects of $L_F(C)$ are lifts of $C$ along $F$.
\item[1.] Morphisms from $(\ovl{C},c: C \cong F(\ovl{C}))$ to $(\ovl{C}',c': C
\cong F(\ovl{C}'))$ are elements of $L_F(1_C: C \lra C|c,c')$ \emph{which are
isomorphisms in $\ovl{\ccc}$}.
\end{enumerate}
\end{enumerate}
\end{definition}

Next we define some natural groupoids for lifting complexes and
cochain maps ``up to homotopy'' relative to fixed
graded lifts.
Let $F: \ovl{\ccc} \lra \ccc$ be an additive functor between linear
categories. There are induced functors $G(F): G(\ovl{\ccc}) \lra G(\ccc)$
between the graded categories, $P(F):P(\ovl{\ccc}) \lra P(\ccc)$ between the
categories of pre-complexes and
$C(F): C(\ovl{\ccc}) \lra C(\ccc)$ between the categories of cochain
complexes. Lifts along $G(F)$ will also be called \emph{graded lifts} whereas
lifts along $C(F)$ will be called \emph{lifts}.

\begin{definition}\label{deflift2}
\begin{enumerate}
\item Consider pre-complexes $(C^{\cdot},d_C), (D^{\cdot},d_D) \in
P(\ccc)$, graded maps $f,g \in \Hom^n(C^{\cdot},D^{\cdot})$ and a homotopy
$H: f \lra g$. Suppose we have lifts $(\ovl{C}^{\cdot},\ovl{d}_C),
(\ovl{D}^{\cdot},\ovl{d}_D)
\in P(\ovl{\ccc})$ along $P(F)$ and graded lifts $\ovl{f}, \ovl{g}$. A
\emph{graded lift of $H$ (along $F$) relative to $\ovl{d}_C,\ovl{d}_D,\ovl{f},
\ovl{g}$} is a graded lift $\ovl{H}$ of $H$ which is a homotopy $\ovl{H}:
\ovl{f} \lra
\ovl{g}$. We consider the following groupoid
$$L_F(H\,|\,\ovl{d}_C,\ovl{d}_D,\ovl{f}, \ovl{g}):$$
\begin{enumerate}
\item[0.] Objects are graded lifts
of $H$ relative to $\ovl{d}_C,\ovl{d}_D,\ovl{f}, \ovl{g}$.
\item[1.] Morphisms from $\ovl{H}$ to $\ovl{H}'$ are graded lifts
$\ovl{0}:\ovl{H} \lra \ovl{H}'$ of $0: H \lra H$ relative to
$\ovl{d}_C,\ovl{d}_D,\ovl{H}, \ovl{H}'$.
\end{enumerate}
\item Consider cochain complexes $(C^{\cdot},d_C), (D^{\cdot},d_D) \in
C(\ccc)$ and a cochain map $f: (C^{\cdot},d_C) \lra (D^{\cdot},d_D)$. For
lifts $(\ovl{C}^{\cdot},\ovl{d}_C),
(\ovl{D}^{\cdot},\ovl{d}_D)
\in C(\ovl{\ccc})$,
we put
$$L_F(f\,|\,\ovl{d}_C,\ovl{d}_D) = L_F(f\,|\,\ovl{d}_C,\ovl{d}_D,0,0),
\text{i.e.}$$
\begin{enumerate}
\item[0.] Objects are lifts of $f$ along $C(F)$ relative to
$(\ovl{C}^{\cdot},\ovl{d}_C), (\ovl{D}^{\cdot},\ovl{d}_D)$.
\item[1.] Morphisms from $\ovl{f}$ to $\ovl{f}'$ are graded lifts of $0: f \lra
f$ relative to $\ovl{d}_C,\ovl{d}_D,\ovl{f}, \ovl{f}'$.
\end{enumerate}
\item Consider a cochain complex
$(C^{\cdot},d)
\in C(\ccc)$ and a graded lift $\ovl{C}^{\cdot}$ of
$C^{\cdot}$.
\begin{enumerate}
\item A \emph{graded lift
of
$d$ along $F$ relative to $\ovl{C}^{\cdot}$} is an element of $L_{G(F)}(d:
C^{\cdot}
\lra C^{\cdot}[1]\,|\,\ovl{C}^{\cdot},\ovl{C}^{\cdot}[1])$.
\item A \emph{lift
of
$d$ along $F$ relative to $\ovl{C}^{\cdot}$} is a graded lift $\ovl{d}$ of $d$
with $\ovl{d}^2 = 0$, i.e. a differential $\ovl{d}$ on $\ovl{C}^{\cdot}$
making $(\ovl{C}^{\cdot}, \ovl{d})$ into a lift of $(C^{\cdot},d)$ along
$C(F)$.
\end{enumerate}
We consider the following bigroupoid $$L_F(d\,| \,\ovl{C}^{\cdot}):$$
\begin{enumerate}
\item[0.] Objects (0-cells) are lifts of $d$ relative to $\ovl{C}^{\cdot}$.
\item[1.] Morphisms (1-cells) from $\ovl{d}$ to $\ovl{d}'$ are lifts of $1:
(C^{\cdot},d) \lra (C^{\cdot},d)$ relative to
$(\ovl{C}^{\cdot},\ovl{d}),(\ovl{C}^{\cdot},\ovl{d}')$ \emph{which are
isomorphisms in $C(\ovl{\ccc})$}.
\item[2.] 2-cells from $\ovl{1}$ to $\ovl{1}'$ are graded lifts of $0: 1 \lra
1$ relative to $\ovl{d}, \ovl{d}', \ovl{1}, \ovl{1}'$.
\end{enumerate}
\end{enumerate}
\end{definition}

\subsection{Obstruction theory}\label{parobs} In this section we give an
obstruction theory for the lift groupoids defined in the previous section
under certain assumptions on $F$. For the additive functor $F: \ovl{\ccc} \lra
\ccc$, let $\Kern (F)$ (resp. $\Kern(F)^2$) be the category-without-identities
with the same objects as
$\ovl{\ccc}$ and containing precisely the $\ovl{\ccc}$-morphisms $f$ with
$F({f}) = 0$ (resp. the compositions of two such morphisms). From now on we
will assume that \emph{$F$ is full} and
\begin{equation}\label{kernil}
\Kern(F)^2 = 0
\end{equation}

This has the following important consequence, which generalizes the well known
fact for rings:

\begin{proposition}\label{isoliftiso}
Suppose $f: C \lra D$ and $g: D \lra C$ are inverse isomorphisms in $\ccc$ and
consider lifts $\ovl{C}$, $\ovl{D}$ of $C$ and $D$ respectively. For every
lift $\ovl{f}: \ovl{C} \lra \ovl{D}$ of $f$, there exists a lift $\ovl{g}:
\ovl{D} \lra \ovl{C}$ of $g$ such that $\ovl{f}$ and $\ovl{g}$ are inverse
isomorphisms. In particular, $\ovl{C}$ and $\ovl{D}$ are isomorphic.
\end{proposition}

\begin{proof}
Consider an arbitrary lift $\ovl{g}'$ of $g$ and suppose $\ovl{f}\ovl{g} - 1 =
\epsilon \in \Kern(F)$. It suffices to change $\ovl{g}'$ into $\ovl{g} =
\ovl{g}'(1 - \epsilon)$.
\end{proof}

In particular, the requirements in Definitions \ref{deflift1}(3) and
\ref{deflift2} that morphisms in the lift groupoids are isomorphisms in
$\ovl{\ccc}$ and $C(\ovl{\ccc})$ respectively are automatically fulfilled.
All (graded) lifts will be along $F$, so we will no longer explicitely say so.
For legibility, we will suppress $F$ in all our notations.

\begin{remark}\label{remark}
If we are only interested in lifting complexes of objects in a certain
subcategory $\ccc' \subset \ccc$, by restricting the codomain of $F$, it
suffices to require that $F$ is full on the closure of $\ccc'$ under
isomorphic objects.
\end{remark}

Consider pre-complexes $(C^{\cdot}, d_C)$ and $(D^{\cdot},d_D)$ with fixed
graded lifts
$\ovl{C}^{\cdot}$ and $\ovl{D}^{\cdot}$.
We define the pre-complex
$$(\mathbf{C}^{\cdot},\ovl{\delta})
= (\mathbf{C}^{\cdot},\ovl{\delta})_{{d}_C,{d}_D}$$
to be the kernel in the
exact sequence of pre-complexes
\begin{equation}\label{defboldc}
\xymatrix{0 \ar[r] &{(\mathbf{C}^{\cdot},\ovl{\delta})} \ar[r] &
{(\Hom^{\cdot}(\ovl{C}^{\cdot},
\ovl{D}^{\cdot}),\ovl{\delta})} \ar[r] &
{(\Hom^{\cdot}(C^{\cdot},D^{\cdot}),\delta)}
\ar[r] &0}
\end{equation}
where $\delta = \delta_{d_C,d_D}$ and $\ovl{\delta} =
\delta_{\ovl{d}_C,\ovl{d}_D}$ for arbitrary graded lifts $\ovl{d}_C,\ovl{d}_D$
of $d_C, d_D$ respectively (i.e. $\ovl{\delta}^n(f) = \ovl{d}_D{f}
- (-1)^n{f}\ovl{d}_C$, see also \S\ref{not}).

\begin{proposition}\label{boldc}
$(\mathbf{C}^{\cdot},\ovl{\delta})$ is a cochain complex which is independent
of the choice of $\ovl{d}_C,\ovl{d}_D$.
\end{proposition}
\begin{proof}
For ${f} \in \mathbf{C}^{\cdot}$, we have
%$\ovl{\delta}^{n+1}(\ovl{\delta}^n({f})) = \ovl{\delta}^{n+1}({f}\ovl{d}_C -
%(-1)^n\ovl{d}_D{f}) = f\ovl{d}_C^2 -(-1)^n\ovl{d}_D{f}\ovl{d}_C
%-(-1)^{n+1}\ovl{d}_D{f}\ovl{d}_C - \ovl{d}_C^2f = f\ovl{d}_C^2 - 
%\ovl{d}_D^2f$.
$\ovl{\delta}\ovl{\delta}(f) = \ovl{d}_D^2f - f\ovl{d}_C^2$.
Since $f$, $\ovl{d}_C^2$ and $\ovl{d}_D^2$ belong to $\Kern(G(F))$, the
expression equals zero by (\ref{kernil}).
Next, other graded lifts of $d_C$ and $d_D$ can be written as $\ovl{d}_C +
\partial_C$ and $\ovl{d}_D +
\partial_D$ for $\partial_C, \partial_D \in \Kern(G(F))$. Hence, for $f \in
\mathbf{C}^{\cdot}$, we have $\delta_{\ovl{d}_C +
\partial_C, \ovl{d}_D +
\partial_D}(f) = \delta_{\ovl{d}_C,\ovl{d}_D}(f) + \partial_D{f}
-(-1)^n{f}\partial_C$ and the last two terms equal zero by (\ref{kernil}).
\end{proof}

The following theorem gives the obstruction theory for lifting homotopies. It
has the obstruction theory for lifting cochain maps as an immediate corollary
(Corollary \ref{obscoc}). For any category $\ddd$, the 
\emph{skeleton} $\Sk(\ddd)$ of $\ddd$ is the
class of all isomorphism classes of $\ddd$-objects.

\begin{theorem}\label{groot}
Consider the following data in $\ccc$:
\begin{itemize}
\item Pre-complexes
$(C^{\cdot},d_C)$ and $(D^{\cdot},d_D)$.
\item Graded maps $f,g: C^{\cdot}
\lra D^{\cdot}$ of degree $n$.
\item A homotopy
$H: f
\lra g$.
\end{itemize}
Suppose we have fixed lifts
$(\ovl{C}^{\cdot},\ovl{d}_C)$ and
$(\ovl{D}^{\cdot},\ovl{d}_D)$ along $P(F)$ of $(C^{\cdot},d_C)$ and
$(D^{\cdot},d_D)$ respectively.
On $\Hom^{\cdot}(\ovl{C}^{\cdot},\ovl{D}^{\cdot})$, put $\ovl{\delta} =
\delta_{\ovl{d}_C,\ovl{d}_D}$. Put $$\mathbf{C}^{\cdot} =
(\mathbf{C}^{\cdot},\ovl{\delta})_{{d}_C,{d}_D}.$$
Suppose we have graded lifts $\ovl{f}, \ovl{g}, \ovl{H}$ of
$f, g, H$ respectively with
$\ovl{\delta}(\ovl{f}) = \ovl{\delta}(\ovl{g})$. Put $L(H) =
L(H\,|\,\ovl{d}_C,\ovl{d}_D,\ovl{f},\ovl{g})$.
\begin{enumerate}
\item
There is an obstruction
$$o_n(H) = o_n(H\,|\,\ovl{d}_C,\ovl{d}_D,\ovl{f},\ovl{g}) = [\ovl{g} - \ovl{f}
- \ovl{\delta}(\ovl{H})]
\in H^n\mathbf{C}$$
with $$o_n(H) = 0 \iff
L(H) \neq \varnothing.$$
\item If $o_n(H) = 0$, the map
$$v_{n-1}: |L(H)|^{2} \rightarrow
H^{n-1}\mathbf{C}: (\ovl{H}, \ovl{H}')
\mapsto [\ovl{H}' - \ovl{H}]$$
satisfies $$v_{n-1}(\ovl{H}, \ovl{H}') = 0 \iff [\ovl{H}] = [\ovl{H}'] \in
\Sk(L(H))$$ and induces an
$H^{n-1}\mathbf{C}$-affine structure on
$\Sk(L(H))$.
\end{enumerate}
\end{theorem}

\begin{proof}
Without loss of generality, we assume that $F(\ovl{C}^{\cdot}) =
C^{\cdot}$ and $F(\ovl{D}^{\cdot}) = D^{\cdot}$.

(1) Clearly, $F(\ovl{g} - \ovl{f} -
\ovl{\delta}(\ovl{H})) = g - f - \delta(H) = 0$ and 
$\ovl{\delta}(\ovl{g} - \ovl{f} -
\ovl{\delta}(\ovl{H})) = \ovl{\delta}(\ovl{g}) - \ovl{\delta}(\ovl{f}) = 0$, so
$\ovl{g} - \ovl{f} -
\ovl{\delta}(\ovl{H})$ is in $Z^n(\mathbf{C}^{\cdot})$. Furthermore,
$L(H\,|\,\ovl{d}_C,\ovl{d}_D,\ovl{f},\ovl{g}) \neq \varnothing$ if and only
if there exists a $\gamma \in \mathbf{C}^{n-1}$ such that $\ovl{H} + \gamma$ is
a homotopy $\ovl{f} \lra \ovl{g}$ or in other words $\ovl{g} - \ovl{f} -
\ovl{\delta}(\ovl{H}) = \ovl{\delta}(\gamma)$ which finishes the proof of
(1).

(2) Since $0:\ovl{H} \lra \ovl{H}'$ is a graded lift of $0: H \lra H$, by
part (1) we have
$o_{n-1}(0:H \lra H\, |\,\ovl{d}_C,\ovl{d}_D,\ovl{H}, \ovl{H}') = [\ovl{H}' -
\ovl{H}]$ which proves the first part of (2). Now it is easily seen that
$$a_{n-1}:|L(H)| \times Z^{n-1}\mathbf{C} \lra |L(H)|: (\ovl{H},\gamma)
\longmapsto \ovl{H} + \gamma$$
defines a strictly transitive action
$\tilde{a}_{n-1}: \Sk(L(H)) \times H^{n-1}\mathbf{C}^{\cdot} \lra \Sk(L(H))$
with difference map $\tilde{v}_{n-1}: \Sk(L(H))^2 \lra
H^{n-1}\mathbf{C}^{\cdot}$ induced by $v_{n-1}$.
\end{proof}

\begin{corollary}\label{obscoc}
Consider the following data in $\ccc$:
\begin{itemize}
\item Pre-complexes
$(C^{\cdot},d_C)$ and $(D^{\cdot},d_D)$.
\item A cochain map $f: (C^{\cdot},d_C) \lra (D^{\cdot},d_D)$ of degree $n$.
\end{itemize}
Suppose we have fixed lifts
$(\ovl{C}^{\cdot},\ovl{d}_C)$ and
$(\ovl{D}^{\cdot},\ovl{d}_D)$ along $P(F)$ of $(C^{\cdot},d_C)$ and
$(D^{\cdot},d_D)$ respectively.
Put $$\mathbf{C}^{\cdot} =
(\mathbf{C}^{\cdot},\ovl{\delta})_{{d}_C,{d}_D}$$
and put
$L(f) = L(f\,|\,\ovl{d}_C,\ovl{d}_D)$.
Suppose we have a graded lift
$\ovl{f}$ of $f$.
\begin{enumerate}
\item There
is an obstruction $$o_{n+1}(f) = o_{n+1}(f\,|\,\ovl{d}_C,\ovl{d}_D) =
[\ovl{\delta}(\ovl{f})]
\in H^{n+1}\mathbf{C}$$
with $$o_{n+1}(f) = 0 \iff L(f) \neq
\varnothing.$$
\item If $o_{n+1}(f) = 0$, the map
$$o_n: |L(f)|^2 \longrightarrow H^n\mathbf{C}: (\ovl{f},\ovl{f}') \longmapsto
[\ovl{f}' - \ovl{f}]$$ satisfies
$$o_n(\ovl{f},\ovl{f}') = 0 \iff [\ovl{f}] = [\ovl{f}'] \in
\Sk(L(f))$$ and
induces an $H^n\mathbf{C}$-affine structure on
$\Sk(L(f))$.
\end{enumerate}
\end{corollary}
\begin{proof}
This is a special case of Theorem \ref{groot} since a cochain map $f$ is a
homotopy
$f: 0
\lra 0$ and we can lift both zeros to zero.
\end{proof}

The following theorem gives the obstruction theory for lifting differentials.

\begin{theorem}\label{groot2}
Consider a cochain complex $(C^{\cdot},d)$ in $\ccc$ with a fixed graded lift
$\ovl{C}^{\cdot}$ of $C^{\cdot}$.
Put $$\mathbf{C}^{\cdot} =
(\mathbf{C}^{\cdot},\ovl{\delta})_{{d},{d}}$$
and put $L(d) = L(d\,|\, \ovl{C}^{\cdot})$.
\begin{enumerate} \item There is an obstruction $$o(d) = o(d\,|\,
\ovl{C}^{\cdot}) = [\ovl{d}^2] \in H^2\mathbf{C}$$
with $$o(d) = 0 \iff L(d) \neq \varnothing.$$
\item If $o(d) = 0$, the map
$$v: |L(d)|^2 \lra H^1\mathbf{C}: (\ovl{d},\ovl{d}') \longmapsto
[\ovl{d}' - \ovl{d}]$$
satisfies $$v(\ovl{d},\ovl{d}') = 0 \iff [\ovl{d}] = [\ovl{d}']
\in \Sk(L(d))$$ and
induces an $H^1\mathbf{C}$-affine structure on $\Sk(L(d))$.
\item If $v(\ovl{d},\ovl{d}') = 0$, the map
$$w: L(d)(\ovl{d},\ovl{d}')^2 \lra H^0\mathbf{C}: (\ovl{1},\ovl{1}')
\longmapsto [\ovl{1}' - \ovl{1}]$$ satisfies
$$w(\ovl{1},\ovl{1}') = 0
\iff [\ovl{1}] = [\ovl{1}'] \in
\Sk(L(d)(\ovl{d},\ovl{d}'))$$ and induces an
$H^0\mathbf{C}$-affine structure on
$\Sk(L(d)(\ovl{d},\ovl{d}'))$.
\end{enumerate}
\end{theorem}

\begin{proof}
Without loss of generality, we assume that $F(\ovl{C}^{\cdot}) =
C^{\cdot}$. On $\Hom^{\cdot}(\ovl{C}^{\cdot},\ovl{C}^{\cdot})$, we put
$\ovl{\delta} =
\delta_{\ovl{d},\ovl{d}}$ for a graded lift $\ovl{d}$ of
$d$.

(1) Clearly, $F(\ovl{d}^2) = d^2 = 0$ and $\ovl{\delta}(\ovl{d}^2) =
\ovl{d}^3 - \ovl{d}^3 = 0$, so $\ovl{d}^2$ is in $Z^2(\mathbf{C}^{\cdot})$.
Furthermore, $L(d\,|\,\ovl{C}^{\cdot}) \neq \varnothing$ if and only if there
exists a
$\partial
\in \mathbf{C}^{1}$ such that $\ovl{d} + \partial$ is a differential on
$\ovl{C}^{\cdot}$ or in other words $\ovl{d}^2 + \ovl{\delta}(\partial) = 0$
which finishes the proof of (1).

(2) By Lemma \ref{boldc}, the differentials $\delta_{\ovl{d},\ovl{d}}$ and
$\delta_{\ovl{d},\ovl{d}'}$ coincide on $\mathbf{C}^{\cdot}$. Consequently,
since $1: \ovl{C}^{\cdot} \lra \ovl{C}^{\cdot}$ is a graded lift of $1:
C^{\cdot} \lra C^{\cdot}$, we have
$o_1(1: C^{\cdot} \lra C^{\cdot}\,|\,\ovl{d}, \ovl{d}') = [\ovl{d}' -
\ovl{d}]$ by Corollary \ref{obscoc}(1) which proves the first part of (2).
Now it is easily seen that
$$a:|L(d)| \times Z^{1}\mathbf{C} \lra |L(d)|: (\ovl{d},\partial)
\longmapsto \ovl{d} + \partial$$defines a strictly transitive action 
$\tilde{a}: \Sk(L(d)) \times
H^1(\mathbf{C}) \lra \Sk(L(d))$ with difference map $\tilde{v}: \Sk(L(d))^2
\lra H^1(\mathbf{C})$ induced by $v$.

(3) This follows from Corollary \ref{obscoc}(2).
\end{proof}

\subsection{Crude lifting lemma}\label{parcrude}

The main theorem of this section, Theorem \ref{comp2}, is entirely inspired by
\cite[Theorem 1.2]{huebs},
\cite[Theorem 3]{markl}. However, to apply these results, we would
have to start with a homotopy equivalence in $\ovl{\ccc}$ and perturb it with
respect to the filtration coming from $F: \ovl{\ccc} \lra \ccc$, whereas we
merely start with a ``homotopy equivalence modulo $F$'' in the first place. In
\cite{markl}, Markl shows that the element
$[Hg - gK] \in H^{-1}\Hom(D^{\cdot},C^{\cdot})$ is the obstruction against the
extension of the given homotopy equivalence to a strong one (\cite[Definition
1]{markl}), for which an Ideal Perturbation Lemma \cite[IPL]{markl} holds.  He
also shows that changing $K$ into
$K' = K + f(Hg - gK)$ kills this obstruction \cite[Theorem 13]{markl},
yielding \cite[Theorem 3]{markl} and the Crude Perturbation Lemma.

\begin{theorem}\label{comp2}
Consider the following data in $\ccc$:
\begin{itemize}
\item Cochain complexes $(C^{\cdot},d_C)$ and $(D^{\cdot},d_D)$.
\item Cochain maps $f: C^{\cdot} \lra D^{\cdot}$ and $g: D^{\cdot} \lra
C^{\cdot}$.
\item Homotopies $H: gf \lra 1_C$ and $K: fg \lra 1_D$.
\end{itemize}
Suppose we have fixed graded lifts $\ovl{C}^{\cdot}$ and $\ovl{D}^{\cdot}$ of
$C^{\cdot}$ and $D^{\cdot}$ respectively.
\begin{enumerate}
\item Suppose we have $\ovl{d}_D \in L(d_D\,|\,\ovl{D}^{\cdot})$.
There exists a $\ovl{d}_C \in L(d_C\,|\,\ovl{C}^{\cdot})$ with $o(f\,|
\ovl{d}_C,\ovl{d}_D) = 0$.
\item
Suppose we have $\ovl{d}_D \in L(d_D\,|\,\ovl{D}^{\cdot})$,
$\ovl{d}_C \in L(d_C\,|\,\ovl{C}^{\cdot})$ and $\ovl{f} \in L(f\,|
\ovl{d}_C,\ovl{d}_D)$. If
\begin{equation}\label{mark}
[Hg - gK] = 0 \in H^{-1}\Hom(D^{\cdot},C^{\cdot})
\end{equation}
then there exist
$\ovl{g} \in L(g\,|\,\ovl{d}_D,\ovl{d}_C)$ and homotopies
$\ovl{H}:\ovl{g}\ovl{f} \lra 1$ and $\ovl{K}: \ovl{f}\ovl{g} \lra 1$ lifting
$H$ and $K$.
\item If we change $K$ into $K' = K + f(Hg - gK)$, then (\ref{mark}) holds.
\end{enumerate}
\end{theorem}

\begin{proof}
Consider $\ovl{d}_D \in L(d_D\,|\,\ovl{D}^{\cdot})$ and take graded lifts
$\ovl{d}_C,\ovl{f},\ovl{g},\ovl{H},\ovl{K}$ where appropriate. We will
abusively denote
$\delta_{\ovl{d}_C,\ovl{d}_C}$,
$\delta_{\ovl{d}_C,\ovl{d}_D}$, $\delta_{\ovl{d}_D,\ovl{d}_C}$ and
$\delta_{\ovl{d}_D,\ovl{d}_D}$ by
$\ovl{\delta}$.
We will gradually change $\ovl{d}_C,\ovl{f},\ovl{g},\ovl{H},\ovl{K}$ until the
required properties hold. At any stage of the proof, we put
$\mu_H = 1 - \ovl{g}\ovl{f} - \ovl{\delta}(\ovl{H}) \in \Kern(G(F))$ and
$\mu_K = 1 - \ovl{f}\ovl{g} -
\ovl{\delta}(\ovl{K}) \in \Kern(G(F))$ for the current graded lifts and we also
have
$\ovl{d}_C^2$,
$\ovl{\delta}(\ovl{f})$ and $\ovl{\delta}(\ovl{g})$ in $\Kern(G(F))$.

First, we will show that $o(d_C\,|\,\ovl{C}^{\cdot}) = 0$.
We have $\ovl{\delta}\ovl{\delta}(\ovl{f}) = \ovl{d}_D^2\ovl{f} -
\ovl{f}{\ovl{d}}_C^2 = -\ovl{f}{\ovl{d}}_C^2$ and consequently
$\ovl{d}_C^2 = \ovl{\delta}(\eta)$ for some $\eta$ with $F(\eta) = 0$ by
Lemma \ref{techlem}(2) below. So from now on, we may and will suppose that
$$\ovl{d}_C^2 = 0.$$ Next we will change $\ovl{d}_C$ into $\ovl{d}_C +
\partial_C$ for some
$\partial_C$ with $F(\partial_C) = 0 = \ovl{\delta}(\partial_C)$ in order
to make $0 = o(f) = o(f\,|\,\ovl{d}_C + \partial_C,\ovl{d}_D) =
[\ovl{\delta}(\ovl{f}) - \ovl{f}\partial_C]$. By Lemma \ref{techlem}(1), $o(f)
= 0$ for $\partial_C = \ovl{g}\ovl{\delta}(\ovl{f})$. From now on, we may and
will suppose that
$$\ovl{\delta}(\ovl{f}) = 0.$$
Next we show that this implies $o(g) = o(g\,|\,\ovl{d}_D,\ovl{d}_C) = 0$.
Indeed, since $\ovl{\delta}(\mu_K) = \ovl{\delta}(\ovl{f})\ovl{g} +
\ovl{f}\ovl{\delta}(\ovl{g}) = \ovl{f}\ovl{\delta}(\ovl{g})$, $o(g) = 0$ by
Lemma \ref{techlem}(2). From now on, we may and will
suppose that
$$\ovl{\delta}(\ovl{g}) = 0.$$
Finally we will change $\ovl{g}$ into $\ovl{g} + \gamma$ in order to make
$o(H) = o(H\,|\,\ovl{d}_C,\ovl{d}_C,(\ovl{g} + \gamma)\ovl{f},1) = 0$ and
$o(K) = o(K\,|\,\ovl{d}_D,\ovl{d}_D,\ovl{f}(\ovl{g} + \gamma),1) = 0$.
We have $o(H) = \mu_H - \gamma\ovl{f}$ and $o(K) = \mu_K - \ovl{f}\gamma$.
By Lemma \ref{techlem}(2), $o(H) = 0$ for $\gamma = \mu_H\ovl{g} +
\ovl{\delta}(\epsilon)$ where $F(\epsilon) = 0$, whereas $o(K) = 0$ for $\gamma
= \ovl{g}\mu_K + \ovl{\delta}(\epsilon)$ where $F(\epsilon) = 0$. Hence $o(H) =
0 = o(K)$ for $\gamma = \mu_H\ovl{g}$ provided
\begin{equation}\label{prov}
0 = [\mu_H\ovl{g} - \ovl{g}\mu_K] = [\ovl{\delta}(\ovl{H}\ovl{g} -
\ovl{g}\ovl{K})].
\end{equation}
But by assumption, $\ovl{H}\ovl{g} -
\ovl{g}\ovl{K} = \ovl{\delta}(\ovl{z}) + \zeta$ with $F(\zeta) = 0$, hence
$\ovl{\delta}(\ovl{H}\ovl{g} -
\ovl{g}\ovl{K}) = \ovl{\delta}(\zeta)$ and (\ref{prov}) holds, which finishes
the proof.
\end{proof}

\begin{lemma}\label{techlem}
With the notations of the proof of Theorem \ref{comp2}, let $\xi$ be a graded
map in
$\ovl{\ccc}$ with
$F(\xi) =
\ovl{\delta}(\xi) = 0$.
\begin{enumerate}
\item We have that any of $(1 - \ovl{g}\ovl{f})\xi$, $\xi(1 - \ovl{g}\ovl{f})$,
$(1 - \ovl{f}\ovl{g})\xi$ and $\xi(1 - \ovl{f}\ovl{g})$ equals
$\ovl{\delta}(\epsilon)$ for some $\epsilon$ with $F(\epsilon) = 0$.
\item Suppose either $\ovl{f}\xi$, $\xi\ovl{f}$, $\ovl{g}\xi$ or $\xi\ovl{g}$
equals
$\ovl{\delta}(\epsilon)$ for some $\epsilon$ with $F(\epsilon) = 0$. Then we
have $\xi = \ovl{\delta}(\eta)$ for some $\eta$ with $F(\eta) = 0$.
\end{enumerate}
\end{lemma}

\begin{proof}
(1) We have $(1 - \ovl{g}\ovl{f})\xi = (\mu_H + \ovl{\delta}(\ovl{H})\xi =
\mu_H\xi + \ovl{\delta}(\ovl{H}\xi) + \ovl{H}\ovl{\delta}(\xi)$, in which the
first term equals zero since $F(\mu_H) = 0 = F(\xi)$ and $\Kern(F)^2 = 0$, the
last term equals zero since $\ovl{\delta}(\xi) = 0$, and in the middle term,
$F(\ovl{H}\xi) = F(\ovl{H})F(\xi) = 0$ which proves (1).
(2) Suppose $\ovl{f}\xi = \ovl{\delta}(\epsilon)$ and $F(\epsilon) = 0$. By
(1), $\xi = \ovl{g}\ovl{f}\xi + \ovl{\delta}(\rho)$ with $F(\rho) = 0$ and
$\ovl{g}\ovl{f}\xi = \ovl{g}\ovl{\delta}(\epsilon) =
\ovl{\delta}(\ovl{g}\epsilon) - \ovl{\delta}(\ovl{g})\epsilon$ in which the
second term equals zero since $F(\ovl{\delta}(\ovl{G})) = 0 =
F(\epsilon)$ and in the first term, $F(\ovl{g}\epsilon) = 0$ which proves (2).
\end{proof}

\begin{corollary}[crude lifting lemma]\label{homdef}\label{crude}
Consider the following data in $\ccc$:
\begin{itemize}
\item Cochain complexes $(C^{\cdot},d_C)$ and $(D^{\cdot},d_D)$.
\item Cochain maps $f: C^{\cdot} \lra D^{\cdot}$ and $g: D^{\cdot} \lra
C^{\cdot}$ homotopy inverse to each other.
\end{itemize}
Suppose we have fixed graded lifts $\ovl{C}^{\cdot}$ and $\ovl{D}^{\cdot}$ of
$C^{\cdot}$ and $D^{\cdot}$ respectively.
\begin{enumerate}
\item Suppose we have $\ovl{d}_D \in L(d_D\,|\,\ovl{D}^{\cdot})$.
There exists a $\ovl{d}_C \in L(d_C\,|\,\ovl{C}^{\cdot})$ with $o(f\,|
\ovl{d}_C,\ovl{d}_D) = 0$.
\item Suppose we have $\ovl{d}_D \in L(d_D\,|\,\ovl{D}^{\cdot})$,
$\ovl{d}_C \in L(d_C\,|\,\ovl{C}^{\cdot})$. For every $\ovl{f} \in L(f\,|
\ovl{d}_C,\ovl{d}_D)$, there exists a $\ovl{g} \in L(g\,|
\ovl{d}_D,\ovl{d}_C)$ such that $\ovl{f}$ and $\ovl{g}$ are homotopy inverse
to each other. In particular, $\ovl{f}$ is a homotopy equivalence.
\end{enumerate}
\end{corollary}

\begin{proof}
This follows from Theorem \ref{comp2}.
\end{proof}

The following proposition is a similar result for homotopies, showing in
particular that the obstructions for lifting homotopies are well-defined up to
homotopy.

\begin{proposition}\label{comp}
Consider the following data in $\ccc$:
\begin{itemize}
\item Cochain complexes
$(C^{\cdot},d_C)$ and $(D^{\cdot},d_D)$.
\item Graded maps $f,g: C^{\cdot}
\lra D^{\cdot}$ of degree $n$.
\item Homotopies
$H, K: f
\lra g$.
\item A homotopy $\Lambda: H \lra K$.
\end{itemize}
Suppose we have fixed lifts $(\ovl{C}^{\cdot},\ovl{d}_C)$ and
$(\ovl{D}^{\cdot},\ovl{d}_D)$ of $(C^{\cdot},d_C)$ and $(D^{\cdot},d_D)$
respectively.
On $\Hom^{\cdot}(C^{\cdot},D^{\cdot})$, put $\delta = \delta_{d_C,d_D}$ and on
$\Hom^{\cdot}(\ovl{C}^{\cdot},\ovl{D}^{\cdot})$, put $\ovl{\delta} =
\delta_{\ovl{d}_C,\ovl{d}_D}$.
\begin{enumerate}
\item Suppose there exist graded lifts $\ovl{f}, \ovl{g}$ of $f,g$ with
$\ovl{\delta}(\ovl{f}) =
\ovl{\delta}(\ovl{g})$. We have
$$o_{n-1}(H\,|\,\ovl{d}_C,\ovl{d}_D,\ovl{f}, \ovl{g}) =
o_{n-1}(K\,|\,\ovl{d}_C,\ovl{d}_D,\ovl{f}, \ovl{g})$$
If these obstructions vanish, then for every $\ovl{H} \in
L(H\,|\,\ovl{d}_C,\ovl{d}_D,\ovl{f}, \ovl{g})$, there exists a $\ovl{K} \in
L(K\,|\,\ovl{d}_C,\ovl{d}_D,\ovl{f}, \ovl{g})$ with
$$o_{n-2}(\Lambda\,|\,\ovl{d}_C,\ovl{d}_D,\ovl{H},\ovl{K}) = 0.$$
\item Suppose $f,g$ are cochain maps.
We have $$o_n(f\,|\,\ovl{d}_C,\ovl{d}_D) = o_n(g\,|\,\ovl{d}_C,\ovl{d}_D)$$
If these obstructions vanish, then for every $\ovl{f} \in
L(f\,|\,\ovl{d}_C,\ovl{d}_D)$ there exists a $\ovl{g} \in
L(g\,|\,\ovl{d}_C,\ovl{d}_D)$ with
$$o_{n-1}(H\,|\,\ovl{d}_C,\ovl{d}_D,\ovl{f},\ovl{g}) = 0.$$
\end{enumerate}
\end{proposition}

\begin{proof}
(1) We have $o_{n-1}(H) - o_{n-1}(K) = [\ovl{\delta}(\ovl{K})
- \ovl{\delta}(\ovl{H})] = [\ovl{\delta}(\ovl{\delta}(\ovl{\Lambda}) -
\beta))]$ for some $\beta \in \mathbf{C}^{n-1}$ which proves the first part
of (1). For the second part, it suffices to take an arbitrary $\ovl{K}' \in
L(K\,|\,\ovl{d}_C,\ovl{d}_D,\ovl{f}, \ovl{g})$ and put $\ovl{K} = \ovl{K}' -
\kappa$ for some $\kappa$ representing
$o_{n-2}(\Lambda\,|\,\ovl{d}_C,\ovl{d}_D,\ovl{H},\ovl{K}')$. (2) is a
special case of (1).
\end{proof}

\section{Lifting in the homotopy category}\label{parliftK}
In this section we will use the results of the previous sections to deduce the
obstruction theory for lifting objects and maps along the functor
$$K(F): K(\ovl{\ccc}) \lra K(\ccc)$$
between the homotopy categories for an \emph{essentially surjective},
\emph{full} additive functor
$F:
\ovl{\ccc}
\lra
\ccc$ with $\Kern(F)^2 = 0$.

\begin{theorem}\label{Klift}
Consider a cochain complex $(C^{\cdot},d)$ in $K(\ccc)$. For any graded lift
$\ovl{C}^{\cdot}$ of $C^{\cdot}$, put $$\mathbf{C}^{\cdot} =
(\mathbf{C}^{\cdot},\ovl{\delta})_{d,d}$$
\begin{enumerate}
\item There is an obstruction $$o(C^{\cdot},d) \in H^2\mathbf{C}$$
with $$o(C^{\cdot},d) = 0 \iff L_{K(F)}(C^{\cdot},d) \neq \varnothing.$$
\item If $o(C^{\cdot},d) = 0$, then $\Sk(L_{K(F)}(C^{\cdot},d))$ is affine
over $H^1\mathbf{C}^{\cdot}$.
\end{enumerate}
\end{theorem}

\begin{proof}
This follows from combining Theorem \ref{groot} and Proposition \ref{psi}
below.
\end{proof}

\begin{theorem}\label{Kliftmap}
Consider a cochain map $f: (C^{\cdot},d_C) \lra (D^{\cdot},d_D)$ between
cochain complexes in $\ccc$ and lifts $(\ovl{C}^{\cdot},\ovl{d}_C),
(\ovl{D}^{\cdot},\ovl{d}_D)$ along $K(F)$.
Put $$\mathbf{C}^{\cdot} =
(\mathbf{C}^{\cdot},\ovl{\delta})_{F(\ovl{d}_C),F(\ovl{d}_D)}$$
\begin{enumerate}
\item There is an obstruction
$$o_{K(F)}(f\,|(\ovl{C}^{\cdot},\ovl{d}_C),(\ovl{D}^{\cdot},\ovl{d_D})) \in
H^1\mathbf{C}$$ with
$$o_{K(F)}(f\,|(\ovl{C}^{\cdot},\ovl{d}_C),(\ovl{D}^{\cdot},\ovl{d}_D)) = 0
\iff L_{K(F)}(f\,|(\ovl{C}^{\cdot},\ovl{d}_C),(\ovl{D}^{\cdot},\ovl{d}_D))
\neq \varnothing.$$
\item Suppose $H^{-1}\Hom(C^{\cdot},D^{\cdot}) = 0$. If
$o_{K(F)}(f\,|(\ovl{C}^{\cdot},\ovl{d}_C),(\ovl{D}^{\cdot},\ovl{d}_D)) = 0$,
then
$\Sk(L_{K(F)}(f\,|(\ovl{C}^{\cdot},\ovl{d}_C),(\ovl{D}^{\cdot},\ovl{d}_D)))$
is affine over $H^0\mathbf{C}$.
\end{enumerate}
\end{theorem}
\begin{proof}
This follows from combining Theorem \ref{groot2} and Proposition \ref{phi}
below.
\end{proof}

Consider a cochain complex $(C^{\cdot},d)$ in $K(\ccc)$ and a graded lift
$\ovl{C}^{\cdot}$ of $C^{\cdot}$. Let $\tilde{L}_F(d\,|\,\ovl{C}^{\cdot})$ be
the bigroupoid associated to the groupoid ${L}_F(d\,|\,\ovl{C}^{\cdot})$, i.e.
\begin{enumerate}
\item[0.] Objects are differentials $\ovl{d}$ on $\ovl{C}^{\cdot}$ making
$(\ovl{C}^{\cdot},\ovl{d})$ into a lift of $(C^{\cdot},d)$ along $C(F)$.
\item[1.] The morphisms from $\ovl{d}$ to $\ovl{d}'$ are given by $\Sk(L_F(1:
C^{\cdot} \lra C^{\cdot}\,|\,\ovl{d}, \ovl{d}'))$.
\end{enumerate}

\begin{proposition}\label{psi}
The natural functor $$\Psi: \tilde{L}_F(d\,|\,\ovl{C}^{\cdot}) \lra
L_{K(F)}(C^{\cdot},d): \ovl{d} \longmapsto (\ovl{C}^{\cdot},\ovl{d})$$
is essentially surjective and full. In particular, it induces a bijection
$$\Sk(\tilde{L}_F(d\,|\,\ovl{C}^{\cdot})) \lra
\Sk(L_{K(F)}(C^{\cdot},d)).$$
If $H^{-1}\Hom(C^{\cdot},D^{\cdot}) = 0$, then $\Psi$ is an equivalence.
\end{proposition}

\begin{proof}
Let us prove essential surjectivity first. A lift of $(C^{\cdot},d = d_C)$
along
$K(F)$ consists of a homotopy equivalence $f: (C^{\cdot},d_C) \lra
(D^{\cdot},d_D)$ in $\ccc$ and a cochain complex
$(\ovl{D}^{\cdot},\ovl{d}_D)$ in $\ovl{\ccc}$ with $F(\ovl{D}^{\cdot}) =
D^{\cdot}$ and $F(\ovl{d}_D) = d_D$. By the crude lifting lemma of Corollary
\ref{crude}, there exist a $\ovl{d}_C \in L_F(d\,|\,C^{\cdot})$ and a
homotopy equivalence $\ovl{f} \in L_F(f\,|\, \ovl{d}_C,\ovl{d}_D)$.
Consequently, $(\ovl{D}^{\cdot},\ovl{d}_D) \cong (\ovl{C}^{\cdot},\ovl{d}_C)$
in $L_{K(F)}(C^{\cdot},d)$.
For the remainder of the proposition, we are to consider for
$\ovl{d}, \ovl{d}' \in {L}_F(d\,|\,\ovl{C}^{\cdot})$ the map
$$\Phi: \Sk(L_F(1\,|\ovl{d}, \ovl{d}')) \lra
L_{K(F)}(1\,|\,(\ovl{C}^{\cdot},\ovl{d}),(\ovl{C}^{\cdot},\ovl{d}')),$$
hence the result follows from Proposition \ref{phi} below.
\end{proof}

\begin{corollary}\label{formsmooth}
Consider a cochain complex $(C^{\cdot},d)$. The
natural functor
$$L_{C(F)}(C^{\cdot},d) \lra L_{K(F)}(C^{\cdot},d)$$
is essentially surjective.\qed
\end{corollary}

\begin{corollary}\label{corpsi}
The first part of Proposition \ref{psi} still holds if $\Kern(F)^n$ = 0.
\end{corollary}

\begin{proof}
This follows from Proposition \ref{psi} by induction.
\end{proof}

\begin{proposition}\label{phi}
Consider a cochain map $f: (C^{\cdot},d_C) \lra (D^{\cdot},d_D)$ between
cochain complexes in $\ccc$ and lifts $(\ovl{C}^{\cdot},\ovl{d}_C),
(\ovl{D}^{\cdot},\ovl{d}_D)$ along $C(F)$.
The natural map
$$\Phi:\Sk(L_F(f\,|\,\ovl{d}_C,\ovl{d}_D)) \lra
L_{K(F)}(f\,|\,(\ovl{C}^{\cdot},\ovl{d}_C),(\ovl{D}^{\cdot},\ovl{d}_D))$$
is surjective.
If $H^{-1}\Hom(C^{\cdot},D^{\cdot}) = 0$, then it is a bijection.
\end{proposition}

\begin{proof}
An element in the image of $\Phi$ is given by a $\ovl{g} \in
L_F(g\,|\,\ovl{d}_C,\ovl{d}_D)$ for some $g$ homotopic to $f$. Suppose $H: f
\lra g$ is a homotopy. By Proposition
\ref{comp}(2), there exists some $\ovl{f} \in
L_F(f\,|\,\ovl{d}_C,\ovl{d}_D)$ for which there exists a homotopy $\ovl{H}:
\ovl{f} \lra \ovl{g}$, proving part one. For part two, suppose
$\Phi(\ovl{f}) = \Phi(\ovl{f}')$. So there exists a homotopy $\ovl{H}:
\ovl{f} \lra \ovl{f}'$ lifting some homotopy $H: f \lra f$. Now $\ovl{f} =
\ovl{f}' \in \Sk(L_F(f\,|\,\ovl{d}_C,\ovl{d}_D))$ if there exists a homotopy
lifting $0: f \lra f$. By Proposition \ref{comp}(1), this is the case if
there exists a homotopy $\Lambda: 0 \lra H$, which finishes the proof.
\end{proof}

\section{Application to linear and abelian deformations}\label{parablin}

In this section, we interpret the results of \S\ref{parliftdifcoc} and
\S\ref{parliftK} for linear and abelian deformations.
Consider surjective ringmaps between coherent, commutative rings $$\ovl{R} \lra
R
\lra R_0 = S$$ with $\Kern(\ovl{R} \lra S) = I$, $\Kern(\ovl{R} \lra R) = J$
and
$IJ = 0$. In particular, $J^2 = 0$ and $J$ is an $S$-module.
\subsection{Linear deformations}\label{parlin}

Let $\CCC_0$ be a fixed flat $S$-linear category, i.e. its hom-sets are flat
$S$-modules. A flat
$R$-deformation of
$\CCC_0$ is an ${R}$-linear functor $\CCC \lra \CCC_0$ inducing an
equivalence
$S
\otimes_R
\CCC
\cong \CCC_0$, where $\CCC$ is flat over $R$. Here $S \otimes_R \CCC$ is the
$S$-linear category with the same objects as ${\CCC}$ and $(S \otimes_{{R}}
{\CCC})(C,C') = S \otimes_R {\CCC}(C,C')$.
Consider flat linear deformations
$$\xymatrix{{\ovl{\CCC}} \ar[r]_F & {\CCC} \ar[r]_{(-)_0} & {\CCC_0}}$$
along the given ringmaps.
Since $J^2 =0$, we have $\Kern(F)^2 = 0$ and $F$ is obviously essentially
surjective and full.

\begin{proposition}\label{propo1}
Consider pre-complexes $(C^{\cdot}, d_C)$ and $(D^{\cdot},d_D)$ in $\CCC$ with
graded lifts
$\ovl{C}^{\cdot}$ and $\ovl{D}^{\cdot}$ along $F$.
The complex $\mathbf{C}^{\cdot}
= (\mathbf{C}^{\cdot},\ovl{\delta})_{{d}_C,{d}_D}$ defined in
(\ref{defboldc}) of \S\ref{parobs} is
$$\mathbf{C}^{\cdot} = J \otimes_S \Hom_{\CCC_0}(C_0^{\cdot}, D_0^{\cdot})$$
where $\Hom_{\CCC_0}(C_0^{\cdot}, D_0^{\cdot})$ is endowed with the
differential $\delta_0 = \delta_{(d_C)_0,(d_D)_0}$.
\end{proposition}

\begin{proof}
The complex $(\Hom_{{\CCC}}({C}^{\cdot},{D}^{\cdot}),{\delta})$ is by
definition isomorphic to the complex
$R\otimes_{\ovl{R}}(\Hom_{\ovl{\CCC}}(\ovl{C}^{\cdot},
\ovl{D}^{\cdot}),\ovl{\delta})$,
so since $\ovl{\CCC}$ is a flat $\ovl{R}$ linear category, the kernel in
(\ref{defboldc}) is given by $J 
\otimes_{\ovl{R}}\Hom_{\ovl{\CCC}}(\ovl{C}^{\cdot},
\ovl{D}^{\cdot}) = J \otimes_S \Hom_{\CCC_0}(C_0^{\cdot}, D_0^{\cdot})$ by
change of rings.
\end{proof}

Consequently, all the results of \S\ref{parliftdifcoc} and
\S\ref{parliftK} can be reformulated using this complex. In particular,
Theorems \ref{Klift} and \ref{Kliftmap} yield the
obstruction theory for lifting along $K(F)$:

\begin{theorem}\label{Klift2}
Consider a cochain complex $C^{\cdot}$ in $K(\CCC)$.
Put $$\mathbf{C}^{\cdot} = J \otimes_S \Hom_{\CCC_0}(C_0^{\cdot},
C_0^{\cdot}).$$
\begin{enumerate}
\item There is an obstruction $$o(C^{\cdot}) \in H^2\mathbf{C}$$
with $$o(C^{\cdot}) = 0 \iff L_{K(F)}(C^{\cdot}) \neq \varnothing.$$
\item If $o(C^{\cdot}) = 0$, then $\Sk(L_{K(F)}(C^{\cdot}))$ is affine
over $H^1\mathbf{C}^{\cdot}$.\qed
\end{enumerate}
\end{theorem}

\begin{theorem}\label{Kliftmap2}
Consider a cochain map $f: C^{\cdot} \lra D^{\cdot}$ between
cochain complexes in $\CCC$ and lifts $\ovl{C}^{\cdot}$, $\ovl{D}^{\cdot}$
along $K(F)$. Put
$$\mathbf{C}^{\cdot} = J \otimes_S \Hom_{\CCC_0}(C_0^{\cdot}, D_0^{\cdot}).$$
\begin{enumerate}
\item There is an obstruction
$$o_{K(F)}(f\,|\ovl{C}^{\cdot},\ovl{D}^{\cdot}) \in
H^1\mathbf{C}$$ with
$$o_{K(F)}(f\,|\ovl{C}^{\cdot},\ovl{D}^{\cdot}) = 0
\iff L_{K(F)}(f\,|\ovl{C}^{\cdot},\ovl{D}^{\cdot})
\neq \varnothing.$$
\item Suppose $H^{-1}\Hom_{\CCC}(C^{\cdot},D^{\cdot}) = 0$. If
$o_{K(F)}(f\,|\ovl{C}^{\cdot},\ovl{D}^{\cdot}) = 0$,
then
$$\Sk(L_{K(F)}(f\,|\ovl{C}^{\cdot},\ovl{D}^{\cdot}))$$
is affine over $H^0\mathbf{C}$.\qed
\end{enumerate}
\end{theorem}

\subsection{Abelian deformations}\label{parab} We start with introducing some
terminology for an $R$-linear \emph{abelian} category $\ccc$. An object $C \in
\ccc$ is called \emph{flat} if the right-exact functor $- \otimes_R C:
\mmod(R) \lra
\ccc: R \longmapsto C$ is exact, and dually, $C$ is called \emph{coflat} if the
left-exact functor $\Hom_R(-,C): \mmod(R) \lra \ccc: R \longmapsto C$ is
exact ($\mmod$ denotes the finitely presented modules). The subcategories of
flat and coflat objects are denoted by
$\mathsf{Fl}(\ccc)$ and $\Cof(\ccc)$ respectively. A (selfdual) notion of
\emph{flatness for abelian categories} was defined in
\cite[Def.3.2]{defab}. An abelian category with enough injectives is
\emph{flat} if its injectives are coflat. In general, a small abelian $\ccc$ is
flat if its category of ind-objects, which is a category with enough
injectives, is flat.
\emph{This notion of flatness is different from the one used in \S\ref{parlin}
for linear categories!} However, we have the following connection:
\begin{enumerate}
\item If $\CCC$ is a small, flat $R$-linear category, then $\Mod(\CCC)
=
\Add(\CCC,\Ab)$ is flat as an abelian $R$-linear category.
\item If $\ccc$ is a flat abelian $R$-linear category, then its category of
injectives $\Inj(\ccc)$ is flat as an $R$-linear category.
\end{enumerate}
Now let $\ccc_0$ be a fixed flat abelian
$S$-linear category. A flat abelian $R$-deformation of $\ccc_0$ is an
$R$-linear functor $\ccc_0 \lra \ccc$ inducing an equivalence $\ccc_0 \cong
\ccc_S$, where $\ccc$ is a flat abelian
$R$-linear category.
Here $\ccc_S$ denotes the category of $S$-objects in $\ccc$, i.e. objects $C$
with an
$S$-structure $S \lra \ccc(C,C)$ extending the $R$-structure
\cite[Def.5.2,\S4]{defab}. Consider flat abelian deformations along the
given ringmaps, together with their adjoints:
\begin{equation}\label{mooi}
\xymatrix{{\ovl{R}}\ar[d]\\ {R} \ar[d]\\ {S}}
\hspace{2cm}
\xymatrix{{\ovl{\ccc}}\ar[d]_{R\otimes_{\ovl{R}}-} \\ {\ccc} \ar[d]_{S
\otimes_R-}\\ {\ccc_0}}
\hspace{1cm}
\xymatrix{{\ovl{\ccc}} \\ {\ccc}\ar[u] \\ {\ccc_0}\ar[u]}
\hspace{1cm}
\xymatrix{{\ovl{\ccc}}\ar[d]^{\Hom_{\ovl{R}}(R,-)} \\ {\ccc}
\ar[d]^{\Hom_R(S,-)}\\ {\ccc_0}}
\end{equation}
All our results for $\Hom_{\ovl{R}}(R,-)$ have of course dual results for
$R\otimes_{\ovl{R}}-$.

\begin{proposition}\label{propo2}
\begin{enumerate}
\item $\Kern(\Hom_{\ovl{R}}(R,-))^2 = 0$.
\item
Consider pre-complexes $(C^{\cdot}, d_C)$ and $(D^{\cdot},d_D)$ in
$P(\Cof(\ccc))$ with graded lifts
$\ovl{C}^{\cdot}$ and $\ovl{D}^{\cdot}$ in $G(\Cof(\ovl{\ccc}))$ along $F =
\Hom_{\ovl{R}}(R,-)$. The complex $\mathbf{C}^{\cdot}
= (\mathbf{C}^{\cdot},\ovl{\delta})_{{d}_C,{d}_D}$ defined in
(\ref{defboldc}) of \S\ref{parobs} is
$$\mathbf{C}^{\cdot} =
\Hom_{\ccc_0}(\Hom_R(J,C^{\cdot}),\Hom_R(S,D^{\cdot})).$$
\end{enumerate}
\end{proposition}

\begin{proof}
(1) For $\ovl{C} \in \ovl{\ccc}$, consider the exact sequence $0 \lra
\Hom_{\ovl{R}}(R,\ovl{C}) \lra \ovl{C} \lra J\ovl{C} \lra 0$ where $J\ovl{C}$
is the image of $\ovl{C} \lra \Hom_{\ovl{R}}(J,\ovl{C})$. A map
$f: \ovl{C}^{\cdot} \lra \ovl{D}^{\cdot}$ has $\Hom_{\ovl{R}}(R,f) = 0$ if
and only if $f$ factors as $\ovl{C}^{\cdot} \lra J\ovl{C} \lra
\Hom_{\ovl{R}}(R,\ovl{D}) \lra \ovl{D}$. Clearly, any composition of two
such maps is zero. (2) With the same argument, the kernel of
(\ref{defboldc}) is given by $\Hom_{\ccc}(J\ovl{C}^{\cdot}, D^{\cdot})$,
and under the flatness assumption on $\ovl{C}^{\cdot}$, $J\ovl{C}^{\cdot} =
\Hom_{\ovl{R}}(J,\ovl{C}) = \Hom_R(J, C^{\cdot})$ by change of rings.
\end{proof}

Consequently, restricting the codomain of $\Hom_{\ovl{R}}(R,-)$ as in
Remark \ref{remark}, all the results of
\S\ref{parliftdifcoc} can be reformulated using this complex,  In
\S\ref{cof}, we give the obstruction theory for lifting coflat objects and
maps between them along
$\Hom_{\ovl{R}}(R,-)$, which yields of course an obstuction theory for lifting
along
$G(\Hom_{\ovl{R}}(R,-))$.
Although we do not get a general obstruction theory for lifting along
$K(\Hom_{\ovl{R}}(R,-)): K(\ovl{\ccc}) \lra K(\ccc)$ or along its restriction
to coflat complexes, we \emph{do} get an obstruction theory for its
restriction to complexes of injectives \emph{if the category $\ovl{\ccc}$,
and hence also $\ccc$, has enough injectives}.
The reason is that in this case, as we will show in
\S\ref{parhominj}, Proposition
\ref{nilpinjlift}, the functor $$\Hom_{\ovl{R}}(R,-): \Inj(\ovl{\ccc}) \lra
\Inj(\ccc)$$
is a \emph{linear} deformation, making lifting along
$K(\Hom_{\ovl{R}}(R,-)): K(\Inj(\ovl{\ccc})) \lra
K(\Inj(\ccc))$ related to both \S\ref{parlin} and this section \S\ref{parab}.

\subsection{Lifting in the homotopy category of injectives}\label{parhominj}
In this section we consider flat abelian deformations as in (\ref{mooi}) of
\S\ref{parab} \emph{with enough injectives} and we put
$$F = \Hom_{\ovl{R}}(R,-): \Inj(\ovl{\ccc}) \lra
\Inj(\ccc).$$

\begin{proposition}\label{nilpinjlift}
Let $\ccc_0 \lra \ovl{\ccc}$ be a flat abelian deformation with enough
injectives along $\ovl{R} \lra S$. The functor $\Hom_{\ovl{R}}(S,-):
\Inj(\ovl{\ccc})
\lra
\Inj(\ccc_0)$
is a \emph{linear} deformation.
\end{proposition}

\begin{proof}
For injective objects $\ovl{E}$, $\ovl{F}$ in $\ovl{\ccc}$, it is easily seen
that
\begin{equation}\label{X}
\Hom_{\ovl{\ccc}}(\Hom_{\ovl{R}}(X,\ovl{E}),\ovl{F}) = X
\otimes_{\ovl{R}}
\Hom_{\ovl{\ccc}}(\ovl{E},\ovl{F})
\end{equation}
for any $X \in \mmod(\ovl{R})$ since $\ovl{E}$ is coflat by assumption on
$\ovl{\ccc}$. Applying this to $X = S$, we obtain
\begin{equation}\label{S}
\Hom_{\ovl{\ccc}}(\Hom_{\ovl{R}}(S,\ovl{E}),\Hom_{\ovl{R}}(S,\ovl{F})) = S
\otimes_{\ovl{R}}
\Hom_{\ovl{\ccc}}(\ovl{E},\ovl{F}).
\end{equation}
So it remains to show that for
every
$\ccc_0$-injective object
$E \in \ccc_0$ there exists an injective $\ovl{\ccc}$-object $\overline{E}$
with
$E
\cong
\Hom_{\ovl{R}}(S,\overline{E})$.
Let $E$ be
an injective object of $\ccc_0$. Take a $\ovl{\ccc}$-monomorphism $m:E
\lra \overline{E'}$ to a $\ovl{\ccc}$-injective. We obtain a
$\ccc_0$-monomorphism
$s:E \lra E'= \Hom_{\ovl{R}}(S,\overline{E'})$. Since $E$ is injective in
$\ccc_S$, we find $r: E' \lra E$ with $rs = 1_E$. This gives us an
idempotent $e = sr: E' \lra E'$. Consider the map $\Hom_{\ovl{R}}(S,-):
(\overline{E'},\overline{E'})
\lra (E',E')$. By (\ref{S}), this map has a nilpotent kernel
$I(\overline{E'},\overline{E'})$. It follows that the idempotent $e$ lifts to
an idempotent $\overline{e}$ in $(\overline{E'},\overline{E'})$. This
idempotent $\overline{e}$ splits as
$\overline{e} = \overline{s}\,\overline{r}$ with $\overline{r}:
\overline{E'} \lra \overline{E}$, $\overline{s}:\overline{E} \lra
\overline{E'}$ for a $\ccc$-injective $\overline{E}$. We now obviously find
an isomorphism $E \cong \Hom_{\ovl{R}}(S,\overline{E})$.
\end{proof}

In accordance with (\ref{X}), both Propositions \ref{propo1} and \ref{propo2}
now yield

\begin{proposition}\label{propo3}
Consider pre-complexes of injectives $(C^{\cdot}, d_C)$ and $(D^{\cdot},d_D)$
in $\ccc$ with graded lifts of injectives
$\ovl{C}^{\cdot}$ and $\ovl{D}^{\cdot}$ along $\Hom_{\ovl{R}}(R,-)$.
The complex $\mathbf{C}^{\cdot}
= (\mathbf{C}^{\cdot},\ovl{\delta})_{{d}_C,{d}_D}$ defined in
(\ref{defboldc}) of \S\ref{parobs} is
$$\begin{aligned}
\mathbf{C}^{\cdot} &= J \otimes_S \Hom_{\ccc_0}(\Hom_R(S,C^{\cdot}),
\Hom_R(S,D^{\cdot}))\\ &= \Hom_{\ccc_0}(\Hom_R(J,C^{\cdot}),
\Hom_R(S,D^{\cdot})).
\end{aligned}$$
\end{proposition}

Theorems \ref{Klift2} and \ref{Kliftmap2} yield:

\begin{theorem}\label{Klift3}
Consider a cochain
complex
$C^{\cdot}$ in
$K(\Inj(\ccc))$. Put $$\mathbf{C}^{\cdot} = \Hom_{\ccc_0}(\Hom_R(J,C^{\cdot}),
\Hom_R(S,C^{\cdot})).$$
\begin{enumerate}
\item There is an obstruction $$o(C^{\cdot}) \in H^2\mathbf{C}$$
with $$o(C^{\cdot}) = 0 \iff L_{K(F)}(C^{\cdot}) \neq \varnothing.$$
\item If $o(C^{\cdot}) = 0$, then $\Sk(L_{K(F)}(C^{\cdot}))$ is affine
over $H^1\mathbf{C}^{\cdot}$.\qed
\end{enumerate}
\end{theorem}

\begin{theorem}\label{Kliftmap3}
Consider a cochain map $f: C^{\cdot} \lra D^{\cdot}$ between
cochain complexes in $\Inj(\ccc)$ and lifts $\ovl{C}^{\cdot}$,
$\ovl{D}^{\cdot}$ along $K(F)$. Put
$$\mathbf{C}^{\cdot} = \Hom_{\ccc_0}(\Hom_R(J,C^{\cdot}),
\Hom_R(S,D^{\cdot})).$$
\begin{enumerate}
\item There is an obstruction
$$o_{K(F)}(f\,|\ovl{C}^{\cdot},\ovl{D}^{\cdot}) \in
H^1\mathbf{C}$$ with
$$o_{K(F)}(f\,|\ovl{C}^{\cdot},\ovl{D}^{\cdot}) = 0
\iff L_{K(F)}(f\,|\ovl{C}^{\cdot},\ovl{D}^{\cdot})
\neq \varnothing.$$
\item Suppose $H^{-1}\Hom_{\ccc}(C^{\cdot},D^{\cdot}) = 0$. If
$o_{K(F)}(f\,|\ovl{C}^{\cdot},\ovl{D}^{\cdot}) = 0$,
then
$$\Sk(L_{K(F)}(f\,|\ovl{C}^{\cdot},\ovl{D}^{\cdot}))$$
is affine over $H^0\mathbf{C}$.\qed
\end{enumerate}
\end{theorem}

\section{Derived lifting}\label{parder}

In this section, we use the results of section \S \ref{parhominj} to obtain
obstruction theories for derived lifting along the adjoints of an
abelian deformation (Theorems \ref{Klift4}, \ref{Kliftmap4} and Theorem
\ref{corollary}). This eventually leads to the obstruction theory for coflat
objects (Theorems \ref{Klift5}, \ref{Kliftmap5}).

\subsection{Comparing lift groupoids}\label{clg}
In the sequel, we will often compare lift groupoids as in Definition
\ref{deflift1} of lifts along different functors. We will use
the following technical tool:

\begin{definition}\label{L}
We will say that a diagram of functors
$$\xymatrix{ \ccc \ar[r]^F \ar[d]_H & \ccc' \ar[d]^{H'} \\
\ddd \ar[r]_G & {\ddd'},}$$
satisfies (L) if the following conditions are fulfilled:
\begin{enumerate}
\item The diagram is commutative up to natural isomorphism.
\item $F$ and $G$ are fully faithful.
\item If $H'(C') \cong G(D)$, then there is a $C \in \ccc$ with $C' \cong
F(C)$.
\end{enumerate}
\end{definition}

\begin{proposition}\label{eqliftgdn}
Suppose a diagram as in Definition \ref{L} satisfies (L).
\begin{enumerate}
\item For $D \in
\ddd$ and $D' \cong G(D)$ in $\ddd'$, there is an equivalence of groupoids
$L_H(D) \lra L_{H'}(D')$.
\item
For $f:D_1 \lra D_2$ in $\ddd$, $g \cong G(f)$ in $\ddd'$, $C_1 \in
L_H(D_1)$, $C_1' \cong F(C_1)$, $C_2 \in L_H(D_2)$ and $C_2' \cong F(C_2)$ and
$g \cong G(f)$ in
$\ddd'$, there is a bijection
$L_H(f\,|\,C_1,C_2) \lra L_{H'}(g\,|\,C_1',C_2')$.
\end{enumerate}
\begin{proof}
The proofs of 1 and 2 are similar and easy. For example for 1, it is convenient
to consider the category $\LLL$ with as objects functors $H$ with a
specified object in the codomain of $H$ and maps between $(H:\ccc \lra \ddd,D)$
and $(H':\ccc' \lra \ddd',D')$ given by 4-tupels $(F, G, \eta, f)$ in
which $F$ and $G$ fit into a square as in Definition \ref{L} (but not
necessarily satisfying (L)), $\eta$ is a natural isomorphism $\eta: 
GH \cong H'F$, and $f$ is an isomorphism
$f: D' \cong G(D)$. There is seen to be a lift functor $L: \LLL \lra \Gd$
mapping $(H,D)$ to $L_H(D)$. If the square satisfies (L), $L(F, G, \eta, f)$
is easily seen to be an equivalence.
\end{proof}
\end{proposition}

\subsection{Derived lifting with enough injectives}\label{parderinj}

In this section we consider flat abelian deformations as in (\ref{mooi}) of
\S\ref{parab} \emph{with enough injectives} and we consider the derived
functor
$$F = R\Hom_{\ovl{R}}(R,-):D^+(\ovl{\ccc}) \lra D^+(\ccc).$$
We can now easily deduce the obstruction theory for $F$
from Theorems \ref{Klift3} and \ref{Kliftmap3}.

\begin{theorem}\label{Klift4}
Consider a cochain
complex
$C^{\cdot} \in D^+(\ccc)$. Put $$\mathbf{C}^{\cdot} =
R\Hom_{\ccc_0}(R\Hom_R(J,C^{\cdot}),
R\Hom_R(S,C^{\cdot})).$$
\begin{enumerate}
\item There is an obstruction $$o(C^{\cdot}) \in H^2\mathbf{C}$$
with $$o(C^{\cdot}) = 0 \iff L_{F}(C^{\cdot}) \neq \varnothing.$$
\item If $o(C^{\cdot}) = 0$, then $\Sk(L_{F}(C^{\cdot}))$ is affine
over $H^1\mathbf{C}^{\cdot}$.
\end{enumerate}
\end{theorem}

\begin{theorem}\label{Kliftmap4}
Consider a cochain map $f: C^{\cdot} \lra D^{\cdot}$ between
cochain complexes in $D^+(\ccc)$ and lifts $\ovl{C}^{\cdot}$,
$\ovl{D}^{\cdot}$ along $F$. Put
$$\mathbf{C}^{\cdot} = R\Hom_{\ccc_0}(R\Hom_R(J,C^{\cdot}),
R\Hom_R(S,D^{\cdot})).$$
\begin{enumerate}
\item There is an obstruction
$$o(f\,|\ovl{C}^{\cdot},\ovl{D}^{\cdot}) \in
H^1\mathbf{C}$$ with
$$o(f\,|\ovl{C}^{\cdot},\ovl{D}^{\cdot}) = 0
\iff L_{F}(f\,|\ovl{C}^{\cdot},\ovl{D}^{\cdot})
\neq \varnothing.$$
\item Suppose $\Ext^{-1}_{\ccc}(C^{\cdot},D^{\cdot}) = 0$. If
$o(f\,|\ovl{C}^{\cdot},\ovl{D}^{\cdot}) = 0$,
then
$\Sk(L_{F}(f\,|\ovl{C}^{\cdot},\ovl{D}^{\cdot}))$
is affine over $H^0\mathbf{C}$.
\end{enumerate}
\end{theorem}

We use the following:

\begin{proposition}\label{eqD1}
In the diagram
$$\xymatrix{{D^+(\ovl{\ccc})} \ar[d]_{F} & {K^{+}(\Inj(\ovl{\ccc}))}
\ar[l]_{\cong} \ar[r] \ar[d] & {K(\Inj(\ovl{\ccc}))}
\ar[d]^{K(\Hom_{\ovl{R}}(R,-))}\\ {D^+({\ccc})} & {K^{+}(\Inj({\ccc}))}
\ar[l]^{\cong} \ar[r] & {K(\Inj({\ccc}))}}$$
both squares  satisfy (L). Consequently, for $C^{\cdot} \in D^{+}(\ccc)$,
there is an isomorphism $C^{\cdot} \cong E^{\cdot}$ with $E^{\cdot} \in
K^{+}(\Inj(\ccc))$ and an equivalence
$$L_{K(\Hom_R(S,-))}(E^{\cdot}) \sim L_F(C^{\cdot})$$
and likewise for maps (see Proposition \ref{eqliftgdn}).
\begin{proof}
Obvious.
\end{proof}
\end{proposition}

\begin{proof}[Proof of Theorems \ref{Klift4} and \ref{Kliftmap4}]
Consider $C^{\cdot}, D^{\cdot} \in D^{+}(\ccc)$ and isomorphisms $C^{\cdot}
\cong E^{\cdot}$ and $D^{\cdot} \cong F^{\cdot}$ with $E^{\cdot}, F^{\cdot} \in
K^{+}(\Inj(\ccc))$. We have
$$\Hom_{\ccc_0}(\Hom_{R}(J,E^{\cdot}), \Hom_{R'}(S,F^{\cdot})) =
R\Hom_{\ccc_0}(R\Hom_{R}(J,C^{\cdot}), R\Hom_{R'}(S,D^{\cdot})),$$ hence the
proof follows by combining Theorems \ref{Klift3} and \ref{Kliftmap3} and
Proposition \ref{eqD1}.
\end{proof}

\begin{remark}
If $\ovl{\ccc}$ is a Grothendieck category, one could hope to deduce an
obstruction theory for the unbounded derived categories $D(\ovl{\ccc}) \lra
D(\ccc)$ from the restriction of $K(\Hom_{\ovl{R}}(R,-))$ to homotopically
injective complexes of injectives.  Unfortunately not every lift of
a homotopically injective complex  is homotopically injective. The
canonical counter example is given by 
$\bar{R}={\mathbb Z}/p^2{\mathbb Z}$, $R={\mathbb Z}/p{\mathbb
  Z}$ and the $R$-complex \def\ZZ{\mathbb{Z}}
\[
\cdots\rightarrow
\ZZ/p\ZZ\xrightarrow{0}\ZZ/p\ZZ\xrightarrow{0}\ZZ/p\ZZ\rightarrow
\cdots
\]
which lifts to the non-homotopically injective complex
\[
\cdots\rightarrow
\ZZ/p^2\ZZ\xrightarrow{p}\ZZ/p^2\ZZ\xrightarrow{p}\ZZ/p^2\ZZ\rightarrow
\cdots
\]
In general it is unclear to us if a homotopically
injective complex in $K(\Inj(\ccc))$  always has a homotopically
injective lift to $K(\Inj(\ovl{\ccc}))$. In the bounded below case, this
problem is overcome by the fact that being a complex of injectives is a
property on the graded level.
\end{remark}

\subsection{Lifting complexes of bounded coflat dimension}\label{parbound}
Consider flat abelian deformations as in
(\ref{mooi}) of
\S\ref{parab} of small abelian categories, and consider the associated
deformations of ind-objects, which have enough injectives. We
will discuss some restrictions of the derived functor
$$R\Hom_{\ovl{R}}(R,-): D^{+}(\Ind \ovl{\ccc}) \lra D^{+}(\Ind \ccc)$$
of the previous section, for which (the restrictions of) Theorems \ref{Klift4}
and
\ref{Kliftmap4} still hold (Theorem \ref{corcor}). Since by
enlarging the universe, we may assume that any category is small, the results
of this section hold for arbitrary abelian categories.

In general, $R\Hom_R(-,-)$ and $\Ext^i_R(-,-)$ are defined as derived
functors in the first argument, fixing the second one. However we will use the
following double interpretation in the sequel:
\begin{proposition}\label{bifunct}
If $\ccc$ is a flat $R$-linear abelian category with enough injectives, then
we have a derived bifunctor
$$D^{-}(\mmod(R)) \times D^{+}(\ccc) \lra D^{+}(\ccc):
(M_{\cdot},C^{\cdot})\longmapsto R\Hom_R(M_{\cdot},C^{\cdot}).$$
If $P_{\cdot} \lra M_{\cdot}$ is a bounded above projective resolution of
$M_{\cdot}$ in $\mmod(R)$, and $C^{\cdot} \lra E^{\cdot}$ is a bounded below
injective resolution of $C^{\cdot}$ in $\ccc$, we have
$$R\Hom_R(M_{\cdot},C^{\cdot}) = \Hom_R(P_{\cdot},C^{\cdot}) \cong
\Hom_R(P_{\cdot},E^{\cdot}) \cong
\Hom_R(M_{\cdot},E^{\cdot}).$$
\begin{proof}
This is just the classical proof, since $\ccc$-injectives are coflat and for
projectives $P$ in $\mmod(R)$, $\Hom_R(P,-)$ is exact.
\end{proof}
\end{proposition}
\begin{definition}\label{defcd}
For a complex $C^{\cdot} \in C(\ccc)$, its \emph{coflat dimension} is defined
to be
$$\mathrm{cd}(C^{\cdot}) = \min\{ n \in \N \, |\, \forall M \in
\mmod(R),
\,
\forall |i|>n \, \,
\Ext^i_{R}(M,C^{\cdot}) = 0\}$$
if such an $n$ exists and $\mathrm{cd}(C^{\cdot}) = \infty$ otherwise.
\end{definition}
Note that since
$\Ext^i_R(M,-)$ in $\ccc$ and $\Ind \ccc$ coincide,
$\mathrm{cd}_{\ccc}(C^{\cdot}) = \mathrm{cd}_{\Ind \ccc}(C^{\cdot})$. We
consider the following full subcategories of
$D^+(\ccc)$ and $D^{+}(\Ind
\ccc)$ respectively:
\begin{itemize}
\item $|D^+_{\mathrm{cd}\leq n}(\ccc)| = \{C^{\cdot}\in D^+(\ccc)\,
|\,\mathrm{cd}(C^{\cdot})
\leq n
\} \,\,\,\, (n \in \N)$
\item $|D^+_{\mathrm{fcd}}(\ccc)| = \{C^{\cdot}\in D^+(\ccc)\,
|\,\mathrm{cd}(C^{\cdot}) < \infty \}$
\item $|D^+_{\ccc ,\mathrm{cd}\leq n}(\Ind \ccc)| = \{C^{\cdot}\in
D^+_{\ccc}(\Ind
\ccc)\, |\, \mathrm{cd}(C^{\cdot})
\leq n  \} \,\,\,\, (n \in \N)$
\item $|D^+_{\ccc ,\mathrm{fcd}}(\Ind \ccc)| = \{C^{\cdot}\in D^+_{\ccc}(\Ind
\ccc)\, |\, \mathrm{cd}(C^{\cdot})
< \infty  \}$
\end{itemize}

\begin{proposition}\label{cdeq}
There is a diagram
$$\xymatrix{{D^+_{\mathrm{cd}\leq n}({\ovl{\ccc}})} \ar[r]^-{\cong}
\ar[d]_{R\Hom_{{\ovl{R}}}^{\mathrm{res}}(R,-)} & {D^+_{{\ovl{\ccc}}
,\mathrm{cd}\leq n}(\Ind {\ovl{\ccc}})}
\ar[d]
\ar[r] & {D^{+}(\Ind {\ovl{\ccc}})}\ar[d]^{R\Hom_{{\ovl{R}}}(R,-)}\\
{D^+_{\mathrm{cd}\leq n}(\ccc)} \ar[r]_-{\cong}  &
{D^+_{\ccc ,\mathrm{cd}\leq n}(\Ind \ccc)} \ar[r] &
{D^{+}(\Ind \ccc)}}$$
in which both squares satisfy (L), and a similar diagram with ``$\mathrm{cd}
\leq n$'' replaced by ``$\mathrm{fcd}$'' with the same property. Consequently,
in both cases, for
$C^{\cdot}
\in D^{+}_{\mathrm{cd}\leq n}(\ccc)$, there is an equivalence
$$L_{R\Hom_{\ovl{R}}^{\mathrm{res}}(R,-)}(C^{\cdot}) \lra
L_{R\Hom_{\ovl{R}}(R,-)}(C^{\cdot})$$ and likewise for maps (see Proposition
\ref{eqliftgdn}).
\end{proposition}

\begin{proof}
For $C^{\cdot} \in D^+_{\ovl{\ccc} ,\mathrm{cd}\leq n}(\Ind \ovl{\ccc})$, we
have
$\mathrm{cd}_{\ccc}(R\Hom_{\ovl{R}}(R,C^{\cdot})) \leq n$ by change of rings
for
$R\Hom$. Also, by the equivalence
$D^b(\ovl{\ccc}) \lra D^b_{\ovl{\ccc}}(\Ind \ovl{\ccc})$,
$\Ext^i_{\ovl{R}}(R,C^{\cdot})$ is in
$\ccc$. This yields the middle vertical arrow. The arrow
$R\Hom^{\mathrm{res}}_{\ovl{R}}(R,-)$ is obtained using the two horizontal
equivalences. Next, we show that the right diagram satisfies (L).
There are two
points to be checked. Suppose $C^{\cdot} \in D^+(\Ind \ovl{\ccc})$ is such that
$R\Hom_{\ovl{R}}(R,C^{\cdot})$ is in $D^+_{\ccc ,\mathrm{cd}\leq n}(\Ind
\ccc)$. First, to show that $H^iC^{\cdot} \in \ovl{\ccc}$, we use that $\Ind
\ovl{\ccc}$ is a locally coherent Ab5 category with $\ovl{\ccc}$ as finitely
presented objects. We use the long exact cohomology sequence
$\dots
\lra \Ext^{i-1}_{\ovl{R}}(J,C^{\cdot}) \lra
\Ext^i_{\ovl{R}}(R,C^{\cdot})
\lra
H^iC \lra \Ext^i_{\ovl{R}}(J,C^{\cdot}) \lra \Ext^{i+1}_{\ovl{R}}(R,C^{\cdot})
\lra
\dots$. By the equivalence $D^b(\ccc) \lra D^b_{\ccc}(\Ind \ccc)$,
$\Ext^i_{\ovl{R}}(J,C^{\cdot}) = \Ext^i_S(J,R\Hom_R(S,C^{\cdot}))$ is finitely
presented, hence so is $H^iC$ as an extension of finitely presented objects.
Next, we need to show that $\mathrm{cd}(C^{\cdot}) \leq n$. Writing an
arbitrary $M \in
\mmod(\ovl{R})$ as an extension of modules in $\mmod(R)$, it follows from the
associated long exact Ext sequence that it suffices to prove
$\Ext^i_{\ovl{R}}(M,C^{\cdot}) = 0$ for $|i| > n$ and $M \in \mmod(R)$. But
this follows from the assumption on $R\Hom_{\ovl{R}}(R,C^{\cdot})$ since
$\Ext^i_{\ovl{R}}(M,C^{\cdot}) =
\Ext^i_R(M,R\Hom_{\ovl{R}}(R,C^{\cdot}))$.
\end{proof}

\begin{theorem}\label{corcor}\label{corollary}
Consider $C^{\cdot} \in D^{+}_{\mathrm{cd}\leq n}(\ccc)$.
Theorems \ref{Klift4} and \ref{Kliftmap4} hold for
$$F = R\Hom_{\ovl{R}}^{\mathrm{res}}(R,-):D^+_{\mathrm{cd}\leq n}(\ovl{\ccc})
\lra D^+_{\mathrm{cd}\leq n}(\ccc)$$
and for
$$F = R\Hom_{\ovl{R}}^{\mathrm{res}}(R,-):D^+_{\mathrm{fcd}}(\ovl{\ccc}) \lra
D^+_{\mathrm{fcd}}(\ccc)$$
\end{theorem}

\begin{proof}
This follows from Theorems \ref{Klift4} and \ref{Kliftmap4} and Proposition
\ref{cdeq}.
\end{proof}

\subsection{Lifting coflat objects}\label{cof}
In this section we consider arbitrary flat abelian deformations as in
(\ref{mooi}) of
\S\ref{parab} and we consider
$$F = \Hom_{\ovl{R}}(R,-): \Cof(\ovl{\ccc}) \lra \Cof(\ccc).$$

\begin{theorem}\label{Klift5}
Consider
$C \in \Cof(\ccc)$. Put $$\mathbf{C}^{\cdot} =
R\Hom_{\ccc_0}(R\Hom_R(J,C),
R\Hom_R(S,C)).$$
\begin{enumerate}
\item There is an obstruction $$o(C) \in H^2\mathbf{C}$$
with $$o(C) = 0 \iff L_{F}(C) \neq \varnothing.$$
\item If $o(C) = 0$, then $\Sk(L_{F}(C))$ is affine
over $H^1\mathbf{C}$.
\end{enumerate}
\end{theorem}

\begin{theorem}\label{Kliftmap5}
Consider a map $f: C \lra D$ in $\Cof(\ccc)$ and lifts $\ovl{C}$,
$\ovl{D}$ along $F$. Put
$$\mathbf{C}^{\cdot} = R\Hom_{\ccc_0}(R\Hom_R(J,C),
R\Hom_R(S,D)).$$
\begin{enumerate}
\item There is an obstruction
$$o(f\,|\ovl{C},\ovl{D}) \in
H^1\mathbf{C}$$ with
$$o(f\,|\ovl{C},\ovl{D}) = 0
\iff L_{F}(f\,|\ovl{C},\ovl{D})
\neq \varnothing.$$
\item If
$o(f\,|\ovl{C},\ovl{D}) = 0$,
then
$\Sk(L_{F}(f\,|\ovl{C},\ovl{D}))$
is affine over $H^0\mathbf{C}$.\qed
\end{enumerate}
\end{theorem}

\begin{proposition}\label{eqcof}
The diagram
$$\xymatrix{{\Cof(\ovl{\ccc})} \ar[r] \ar[d]_{F}&
{\D^+(\Ind \ovl{\ccc})}
\ar[d]^{R\Hom_{\ovl{R}}(R,-)} \\ {\Cof(\ccc)} \ar[r] &
{\D^+(\Ind \ccc)}}$$
satisfies (L).
Consequently, for $C \in \Cof(\ccc)$, there is an equivalence
$$L_F(C) \lra L_{R\Hom_{\ovl{R}}(R,-)}(C)$$
and likewise for maps.
\end{proposition}
\begin{proof}
By Proposition \ref{cdeq} for $n = 0$, it suffices to note that there is a
diagram
$$\xymatrix{{\Cof(\ovl{\ccc})} \ar[r]^{\cong} \ar[d]_{F}&
{\D^+_{\mathrm{cd} \leq 0}(\ovl{\ccc})}
\ar[d]^{R\Hom^{\mathrm{res}}_{\ovl{R}}(R,-)} \\ {\Cof(\ccc)} \ar[r]_{\cong} &
{\D^+_{\mathrm{cd} \leq 0}(\ccc)}}$$
which obviously satisfies (L).
\end{proof}

\begin{proof}[Proof of Theorems \ref{Klift5} and \ref{Kliftmap5}]
This follows by combining Theorems \ref{Klift4} and \ref{Kliftmap4} and
Proposition \ref{eqcof}.
\end{proof}

\section{Appendix: Miniversal deformations}

In this appendix we prove the existence of miniversal deformations
in the classical setting of \cite{schlessinger} using \cite[Theorem
2.11]{schlessinger}. The results in \S \ref{versdefdif}, \S \ref{vershomdef}
are well known \cite{Laudal}. The results in \S \ref{versderdef} can be found
in \cite{Inaba} for the derived category of coherent sheaves over a
projective variety and in \cite{defder} for the derived
category of a profinite group.

Let
$S = R_0 = k$ be a field, let
$\hat{\mathbf{C}}$ be the category of complete noetherian local $k$-algebras
$(A,\mathrm{m})$ with residue field $k$ and let $\mathbf{C}$ be its
subcategory of artinian rings. Let $F: \mathbf{C} \lra \Set$ be a functor
such that $F(k)$ is a singleton. Recall that a
\emph{hull} for $F$ \cite[Def 2.7]{schlessinger} is a natural transformation
$\eta: H = \hat{\mathbf{C}}(R,-) \lra F$ (for some
$R \in
\hat{\mathbf{C}}$) such that 
\begin{enumerate}
\item[(H1)] $\eta$ is formally smooth \cite[Def
2.2]{schlessinger}; i.e. every surjective $\mathbf{C}$-map $R' \lra R$
induces a surjection $H(R') \lra H(R) \times_{F(R)} F(R')$.
\item[(H2)] $\eta_{k[\epsilon]/(\epsilon^2)}: H(k[\epsilon]/(\epsilon^2))
\lra F(k[\epsilon]/\epsilon^2)$ is a bijection.
\end{enumerate}
If we extend $F$ to
$\hat{\mathbf{C}}$ by putting $\hat{F}((A,\mathrm{m})) = \mathrm{proj} \lim
F(A/\mathrm{m}^n)$, a hull for $F$ corresponds to an element $\zeta = \eta(1)
\in \hat{F}(R)$, which is called a \emph{miniversal deformation} of the
unique element of $F(k)$. If $\eta$ is a natural isomorphism, $F$ is called
\emph{pro-representable}, and  in this case $\zeta$ is a \emph{universal
deformation}. Schlessingers conditions \cite[Theorem 2.11]{schlessinger} for
the existence of a hull are
\begin{enumerate}
\item[(S1)] If $R' \lra R$ is a surjective $\mathbf{C}$-map with kernel of
dimension $1$, and if
$R''
\lra R$ is any $\mathbf{C}$-map, then the map
\begin{equation}\label{map}
F(R' \times_R R'') \lra F(R') \times_{F(R)} F(R'')
\end{equation}
is surjective.
\item[(S2)] The map (\ref{map}) is bijective when $R' \lra R$ is
$k[\epsilon]/\epsilon^2 \lra k$.
\item[(S3)] The tangent space $F(k[\epsilon]/(\epsilon^2))$ is a finite
dimensional $k$-vector space.
\end{enumerate}
If in addition the maps in (S1) are bijective, then $F$ is pro-representable.

\subsection{Deformations of differentials}\label{versdefdif}

Let $\CCC$ be a fixed $k$-linear category. For $R \in \mathbf{C}$, we
consider the trivial $R$-deformation $F_R: R \otimes_k \CCC \lra \CCC$ of
$\CCC$. For a fixed complex $C^{\cdot} = (C^{\cdot},d) \in
C(\CCC)$, we consider $C^{\cdot}_R = C^{\cdot}$ as fixed \emph{graded} lift
to $R \otimes_k \CCC$. Put $L(R) = L_{F_R}(d\,|\,
C^{\cdot})$ as in Definition \ref{deflift2}(3).
Consequently,
$|L(R)|$ contains all lifts of $d$ to a differential $d_R$ on $C^{\cdot}$ in
$R
\otimes_k \CCC$, and in $\Sk(L(R))$ two such lifts $d_R$ and $d_R'$ are
equivalent if there exists a lift of $1: d \lra d$ to an isomorphism $d_R
\lra d_R'$. We will call the objects of $L(R)$ \emph{($R$-)deformations of
$d$}. We consider the functor
$$
F: \mathbf{C} \lra \Set: R \longmapsto \Sk(L(R)).
$$

\begin{proposition}\label{smooth}
If $\mathrm{dim}_k(K(\CCC)(C^{\cdot},C^{\cdot}[1])) < \infty$,
then $F$ has a hull; in other words, the differential $d$ has a miniversal
deformation.
\end{proposition}

\begin{proof}
Consider the commutative diagram
$$\xymatrix{{|L(R' \times_R R'')|} \ar[d]_{\alpha} \ar[r]^-{\gamma} &
{|L(R')|
\times_{|L(R)|} |L(R'')|} \ar[d]^{\beta} \\ {\Sk(L(R' \times_R R''))}
\ar[r]_-{\delta} & {\Sk(L(R'))
\times_{\Sk(L(R))} \Sk(L(R''))}}$$
By Lemma \ref{cond}(1), $\gamma$ is a bijection and it easily follows from
Lemma \ref{cond}(2) that $\beta$ is surjective if $R' \lra R$ is surjective
with kernel of dimension $1$. Consequently, (H1) holds for $F$. For $R = k$,
both
$|L(R)|$ and
$\Sk(L(R))$ are singletons, which easily implies that $\delta$ is a
bijection hence (H2) holds for
$F$. Finally, $(H3)$ follows from Proposition \ref{propo1}, Theorem
\ref{groot2} and the assumption.
\end{proof}

We have used the following Lemma.

\begin{lemma}\label{cond}
Consider the functor $$F_0: \mathbf{C} \lra \Set: R \longmapsto |L(R)|$$
and the canonical natural transformation $$\mu: F_0 \lra F.$$
\begin{enumerate}
\item For arbitrary $\mathbf{C}$-maps $R' \lra R$ and $R'' \lra R$ , the map
$$
F_0(R' \times_R R'') \lra F_0(R') \times_{F_0(R)} F_0(R'')
$$
of (S1) is bijective.
\item $\mu$ is formally smooth.
\end{enumerate}
\end{lemma}

\begin{proof}
(1) By flatness of $\CCC(C,D)$ over $k$, the canonical
$$(R' \times_R R'') \otimes_k \CCC(C,D) \lra R' \otimes_k \CCC(C,D)
\times_{R \otimes_k \CCC(C,D)}R'' \otimes_k \CCC(C,D)$$
is an isomorphism of $k$-modules. Endowing the right hand side with
componentwise compositions, there results an isomorphism of categories
\begin{equation}\label{iso}
(R' \times_R R'') \otimes_k \CCC \lra R' \otimes_k \CCC
\times_{R \otimes_k \CCC}R'' \otimes_k \CCC.
\end{equation}
Consequently, under (\ref{iso}), a graded lift $d'''$ of $d$ to $(R'
\times_R R'')
\otimes_k
\CCC$ corresponds to a couple $(d',d'')$ of graded lifts $d'$ of $d$ to
$R' \otimes_k \CCC$ and $d''$ of $d$ to $R'' \otimes_k \CCC$ with the
same image in $R \otimes_k \CCC$. Furthermore, ${d'''}^2$ corresponds to
$({d'}^2,{d''}^2)$, which finishes the proof of (1).
(2) We are to show that a surjection $R \lra S$ in
$\mathbf{C}$ induces a surjection $F(R) \lra F(S) \otimes_{F_0(S)} F_0(R)$.
So suppose we have lifted differentials $d_S$ on $C^{\cdot}_S$, $d_R$ on
$C^{\cdot}_R$ with an isomorphism $f: d_S \lra (d_R)|_S$ lifting $1$. Then
there exists a differential $d_R'$ lifting $d_S$ and an isomorphism $d_R'
\lra d_R$ lifting $f$ (this follows by induction on the kernel of $R \lra
S$ from Corollary \ref{crude} and Proposition \ref{isoliftiso}).
\end{proof}

\subsection{Homotopy deformations}\label{vershomdef}
Let $\CCC$ and $F_R: R \otimes_k \CCC \lra \CCC$ be as in the previous
section and consider $K(F_R): K(R \otimes_k \CCC) \lra K(\CCC)$. For
$C^{\cdot} \in K(\CCC)$, we will call the objects of $L_{K(F_R)}(C^{\cdot})$
\emph{homotopy ($R$-)deformations of $C^{\cdot}$}.
We consider the functor 
$$F_1: \mathbf{C} \lra \Set: R \longmapsto \Sk(L_{K(F_R)}(C^{\cdot})).$$

\begin{proposition}\label{F1} There is a natural isomorphism of functors $F_0
\cong F_1$. Consequently, if
$\mathrm{dim}_k(K(\CCC)(C^{\cdot},C^{\cdot}[1])) <
\infty,$  then $F_1$ has a hull; in other words, $C^{\cdot}$ has a
miniversal homotopy deformation.
\end{proposition}

\begin{proof}
For $R \in \mathbf{C}$, $R \lra k$ has a nilpotent kernel hence the natural
bijections $\Sk(L_{F_R}(d_0\,|\,
C^{\cdot}_0)) \lra \Sk(L_{K(F_R)}(C^{\cdot}_0))$ follow by induction from
Proposition \ref{psi}. The remainder of the statement follows from
Proposition \ref{smooth}.
\end{proof}

\subsection{Derived deformations}\label{versderdef} Let $\ccc$ be a fixed
small flat abelian
$k$-linear category. For $R \in \mathbf{C}$, we consider the trivial abelian
deformation
$\ccc \lra \ccc_R: C \longmapsto (C,R \lra k \lra \ccc(C,C))$ of
$\ccc$ and its right adjoint $\Hom_R(k,-):\ccc_R \lra
(\ccc_R)_k \cong \ccc$.
We consider
$$R\Hom_R^{\mathrm{res}}(k,-): D^+_{\mathrm{fcd}}(\ccc_R) \lra
D^+_{\mathrm{fcd}}(\ccc)$$
as in \S \ref{parbound}. For a fixed $C^{\cdot} \in
D^+_{\mathrm{fcd}}(\ccc)$, we will call the objects of
$L_{R\Hom_R^{\mathrm{res}}(k,-)}(C^{\cdot})$ \emph{derived ($R$-)deformations
of $C^{\cdot}$}. We consider
$$F_2: \mathbf{C} \lra \Set: R \longmapsto
\Sk(L_{R\Hom_R^{\mathrm{res}}(k,-)}(C^{\cdot})).$$

\begin{proposition}\label{11}
If $\mathrm{dim}_k(\Ext_{\ccc}^1(C^{\cdot},C^{\cdot})) \leq \infty$, then
$F_2$ has a hull; in other words, $C^{\cdot}$ has a miniversal derived
deformation. 
\end{proposition}

\begin{remark}\label{222}
In exactly the same way, the functor describing derived deformations of
\emph{bounded} coflat dimension and the functor describing coflat
deformations of objects have a hull. The latter is shown in 
\cite[Prop. E1.11]{AZ2} (for noetherian objects).
\end{remark}

%Next we consider
%$$\Hom_R(k,-): \Cof((\ccc)_R) \lra \Cof(\ccc)$$
%as in \S \ref{cof}, and for a fixed $C \in
%\Cof(\ccc)$ we consider
%$$F_3: \mathbf{C} \lra \Set: R \longmapsto
%\Sk(L_{\Hom_R(k,-)}(C)).$$
%
%\begin{proposition}\label{22}
%If $\mathrm{dim}_k(\Ext_{\ccc}^1(C,C) = 0$, then
%$F_3$ is pro-representable.
%\end{proposition}

\begin{proof}
Consider $\Hom_R(k,-): \Inj(\Ind\ccc_R) \lra \Inj(\Ind\ccc)$
and let $E^{\cdot}$ be an injective resolution in
$\Ind\ccc$ of $C^{\cdot}$. By induction on the
nilpotent kernel of $R \lra k$, Propositions
\ref{eqD1},
\ref{cdeq} yield a natural
isomorphism $F_3 \cong F_2$ for
$$F_3: \mathbf{C} \lra \Set: R \longmapsto
\Sk(L_{K(\Hom_R(k,-))}(E^{\cdot})).$$
By Lemma \ref{injinj} below, we have $\Inj(\Ind\ccc_R)
\cong
\Inj((\Ind\ccc)_R) \cong R \otimes_k \Inj(\Ind\ccc)$ hence the
results follows from Proposition \ref{F1}. 
\end{proof}

\begin{lemma}\label{injinj}
For an abelian $k$-linear category $\ccc$ with enough injectives, the functor
$\Hom_k(R,-): \ccc \lra \ccc_R$ induces an isomorphism
$$R \otimes_k \Inj(\ccc) \cong \Inj(\ccc_R).$$
\end{lemma}

\begin{proof}
The composition of forgetful functors $\ccc \cong (\ccc_R)_k \lra
\ccc_R \lra \ccc$ is (naturally isomorphic to) the identity, hence
the same holds for the composition of right adjoints
$$\xymatrix{{\ccc} \ar[d]^{\Hom_k(R,-)} \\ {\ccc_R}
\ar[d]^{\Hom_R(k,-)} \\ {\ccc \cong (\ccc_R)_k}}$$
Consequently, for an injective $E \in \ccc$, $\Hom_k(R,E) \in \ccc_R$
is the unique (up to isomorphism) lift of $E$ along $\Hom_R(k,-)$ (see
Theorem
\ref{Klift5}). For injectives $E, F \in \ccc$, we have
$$R \otimes _k \ccc(E,F) \cong \ccc(\Hom_k(R,E),F) \cong
\ccc_R(\Hom_k(R,E),\Hom_k(R,E))$$
whence the result follows.
\end{proof}

%\bibliography{Liftcrit}
%\bibliographystyle{amsabbrv}
\ifx\undefined\bysame
\newcommand{\bysame}{\leavevmode\hbox to3em{\hrulefill}\,}
\fi

\end{document}